\documentclass[a4paper,12pt]{amsart}
\usepackage[utf8]{inputenc}
\usepackage{amssymb} % librerias ams
\usepackage{amsthm}
\usepackage{amsmath}
\usepackage{amsbsy}
\usepackage{graphicx}
\usepackage{subfig}
\usepackage{dsfont}
\usepackage{mathtools}
\usepackage{color}
\definecolor{light}{gray}{0.8}
\newcommand{\bff}{\textbf{f}}

\newcommand{\bs}{\boldsymbol}

\newcommand{\bsk}{{\bs k}}
\newcommand{\bsl}{{\bs l}}
\newcommand{\bsM}{{\bs M}}
\newcommand{\bsS}{{\bs S}}
\newcommand{\bsu}{{\bs u}}
\newcommand{\bsx}{{\bs x}}

\newcommand{\bsmu}{{\bs \mu}}

\newcommand{\bssi}{{\bs \sigma}}
\newcommand{\bsSi}{{\bs \Sigma}}

\newcommand{\bsxi}{{\bs \xi}}
\newcommand{\bsom}{{\bs \omega}}

\newcommand{\E}{\text{E}}
\newcommand{\Cov}{\text{Cov}}
\newcommand{\Var}{\text{Var}}

\def\C{\textrm{C}}
\def\E{\textrm{E}}
\def\P{\textrm{P}}
\def\S{\textrm{S}}
\def\c#1{[\textrm{#1}]}
\def\equil#1#2{%
  \quad
  \mathop{\rightleftarrows}\limits^{k_{#1}}_{k_{#2}}%
  \quad
}
\def\product#1{%
  \quad
  \mathop{\rightarrow}\limits^{k_{#1}}%
  \quad
}

\begin{document}

%\begin{frontmatter}

\title[PCE for general multivariate distributions]{Polynomial Chaos Expansion for general multivariate distributions with correlated variables}
\author{Maria Navarro, Jeroen Witteveen, Joke Blom$\,^\star$}

\date{\today}

\begin{abstract}
Recently, the use of Polynomial Chaos Expansion (PCE) has been increasing to study the uncertainty
in mathematical models for a wide range of applications and several extensions of the 
original PCE technique have been developed to deal with some of its limitations.
But as of to date PCE methods still have the restriction that the random variables have to be statistically 
independent.
This paper presents a method to construct a basis of the probability space of orthogonal polynomials 
for general multivariate distributions with correlations between the random input variables. 
We show that, as for the current PCE methods, the statistics like mean, variance and Sobol' indices can be
obtained at no significant extra postprocessing costs. We study the behavior of the proposed method for a range of correlation
coefficients for an ODE with model parameters that follow a bivariate normal distribution.  In all cases 
the convergence rate of the proposed 
method is analogous to that for the independent case. Finally, we show, for a canonical enzymatic reaction,
how to propagate experimental errors through the process of fitting parameters to a probabilistic
distribution of the quantities of interest, and we demonstrate the significant difference in the results assuming
independence or full correlation compared to taking into account the true correlation.

\vspace{2mm}
%\begin{keyword}

\vspace{2mm}

\noindent CWI, Science Park 123, 1098 XG Amsterdam, the Netherlands

\noindent $^\star$Corresponding author. Tel.:+31 205924263

\noindent {\it Email addresses.} {\scriptsize \tt{maria.navarro@cwi.nl}, \tt{jeroen.witteveen@cwi.nl}, \tt{joke.blom@cwi.nl}}

\end{abstract}

\maketitle
%\end{frontmatter}

\section{Introduction}
To describe real-life phenomena we often make use of deterministic mathematical models. 
Typical constituents of such models are assigned a definite value and we seek a deterministic solution
to the problem.
In reality, however, those phenomena will almost always have uncertain
components: unknown parameters, imprecise experimental data, etc.. A precise mathematical description should
reflect such uncertainties. In other words, the parameters of a mathematical model of a real-life problem 
possess randomness with, most likely, some degree of correlation between them as well.
In studying the propagation  of the uncertainty through the model and its effect on the final solution 
the great majority of the current techniques usually ignores  the correlations between random inputs.

The simplest approach to study the model uncertainty is to apply a Monte Carlo sampling \cite{MC},
where - once a probability density function (pdf) for the random inputs is assumed based on the
a priori knowledge about them - the mean and other characteristics
can be estimated from the output distribution, by sampling repeatedly from the assumed probability density
function and simulating the model for each sample.
Although with this approach it is possible to consider the correlations between variables,
it yields reasonable results only if the number of samples is quite large, requiring a great 
computational effort. Moreover, the order of convergence is merely $\sqrt{N}$, where $N$ is the number of 
samples. To decrease the computational effort, several modifications leading to new methods were introduced:
Latin hypercube sampling \cite{LHS}, the Quasi-Monte Carlo (QMC) method \cite{QMC},
the Markov Chain Monte Carlo method (MCMC) \cite{MCMC}, the Response Surface method (RSM) \cite{RSM}, etc..
%\JGB{[Which one is considered ``state-of-the-art''?]}
Alternatively, deterministic methods to study parameter sensitivity have also been developed, e.g., 
perturbation methods\cite{KleiberHien1992}, local expansion-based methods like Taylor series, and so on.

In this paper we consider Polynomial Chaos Expansion (PCE), which experienced an increasing relevance 
during the last years. The method is based on Wiener's \cite{PCE} homogeneous chaos theory published in 
1938. Cameron and Martin \cite{CandM} proved convergence for the
classical Wiener-Hermite PC expansions based on the Hermite polynomial functionals  in terms
of Gaussian random variables. In the sixties of the last century some successful applications appeared 
\cite{classic1,classic2}. For non-Gaussian random variables, however, PCE was very slow,
leading to a decrease of interest in the method.
To solve this problem, the PCE method was extended 
%in \cite{XK} Xiu and Karniadakis introduced generalized Polynomial Chaos Expansion (gPCE) that generalized the PCE
to polynomials of the Askey scheme in \cite{XK}.
It sparked again new interest for this method, but application tasks demanded further
adaptation of the method to general distributions. 
For this reason, several PCE extensions have been developed in the last years \cite{aPCE,g2,GSPC2008,aPCE2007,G-S}.
A theoretical framework for arbitrary multivariate probability measures was laid out in \cite{g1}.
%as: the multi-element generalized polynomial chaos (ME-gPCE) \cite{ME-gPCE},
% arbitrary polynomial chaos expansion (aPCE) \cite{aPCE}, Gram-Schmidt Polynomial Chaos \cite{G-S}, and
%the methods in \cite{g1,g2} to name just a few.
These polynomial chaos expansion methods have been shown to be effective in a
large amount of applications in different fields as
can be seen in the literature:
air-water flows in soil \cite{water}, chemical reactions \cite{react},
fluid dynamics \cite{fluid,FD}, stability and control \cite{stcon}, 
%fluid mechanics \cite{fluid},
etc.. New and 
promising fields could be (bio)chemistry and artificial intelligence, e.g., uncertain case-based reasoning
\cite{CBRGolosnoy2008,CBRMaria2009}. 
However, the method still has its limitations: {\em (i)} for large numbers of random variables, PCE 
becomes - unacceptably - computationally expensive, and {\em (ii)} the random variables have to be 
statistically independent. 
Therefore, usually the uncertain inputs are assumed to be independent or fully correlated.
To apply PCE when the inputs are linearly correlated, some methods have been 
proposed: linear transformations \cite{Nataf,Rosenblatt}, K-L expansion \cite{KL} or proper
orthogonal decomposition \cite{OrD}. All of these fix the problem by applying transformations to remove the 
correlations, which increases the complexity of the problem and degrades the convergence rate 
because of increased nonlinearity \cite{EldredEtAl2008}.

\medskip
In this paper, we propose a - more fundamental - solution for the latter of those limitations, since
correlations can have strong dynamical effects on the final solution as we will also show. 
We present a method to construct an orthogonal polynomial basis for any general multivariate 
distribution, including those with correlated random variables.

The outline of this paper is as follows: In Section 2 the extension is presented of the
polynomial chaos expansion to general multivariate distributions and the
derivation of the statistics of
the Quantities of Interest (QoI) - like mean, variance, and Sobol' indices - is discussed.
In Section 3 we show for an ODE dependent on a bivariate normal distribution the effect of the strength of 
the correlation - from uncorrelated to fully correlated - on the solution. We also show that in all these cases
the convergence rate for increasing expansion order is the same, i.e., 
the proposed method shows the same favourable convergence rate as the 
original PC method for an uncorrelated multivariate normal distribution.
We then apply the  method on a simple but realistic example from
biochemistry, viz. the canonical enzymatic reaction \cite{BrownMM02}, that has been used in a previous 
paper \cite{FEBS} to display
the effect of noisy data on the reliability of the estimated parameters. We show that {\em true} propagation
of the uncertainty in the parameters - i.e., including the correlation - effects the uncertainty in the QoI,
the concentration of the product, significantly.
Section 4 contains a discussion and concluding remarks.

\section{PCE for multivariate arbitrarily distributed input variables}
In this section we extend the polynomial chaos expansion for arbitrarily distributed
{\em independent} random variables to the 
{\em multivariate} arbitrarily distributed PC expansion.

\subsection{Polynomial chaos expansion}
Consider the following stochastic equation in the probability space
$(\Omega,\mathcal{A},P)$ where
$\Omega$ is the event space, $\mathcal{A}\subseteq 2^{\Omega}$ its $\sigma$-algebra and $P$ its probability measure
\begin{equation}
\mathcal{L}(\bsx,\bsxi(\bsom); \bsu )=\bff(\bsx,\bsxi(\bsom)),\;\;\bsx\in\text{X and }\bsom\in\Omega,
\label{equ:esto}
\end{equation}
where $\bsu(\bsx,\bsxi(\bsom))$ is the stochastic solution vector,
$\mathcal{L}$ is a %(not necessarily linear) 
differential operator, and $\bff$ a %nonlinear 
source function;
$\bsx$ is the vector of {\em deterministic} input variables describing, e.g., time or space, and
$\bsxi(\bsom)$ is the $n$-dimensional vector of {\em random} input variables, with joint probability density
function $\rho(\bsxi(\bsom))$; $\bsxi$ typically contains the uncertainties in model parameters
or initial and boundary conditions.
In order to make the notation less cumbersome, we denote the realization
of a random vector $\bsxi(\bsom)$, for $\bsom \in \Omega$, by $\bsxi \in \Xi$, with $\Xi$ the support of the pdf.

In the polynomial chaos method the random input $\bsxi$ and the solution $\bsu$ are expanded into a 
series of polynomials
\begin{eqnarray}
\bsxi&=&\sum_{i=0}^\infty \bsxi_i \Phi_i(\bsxi)\nonumber\\
\bsu(\bsx,\bsxi) &=& 
\sum_{i=0}^\infty \bsu_i(\bsx)\Phi_i(\bsxi),
\label{equ:expansion}
\end{eqnarray}
separating the {\em deterministic} and the {\em random} variables. The
$\Phi_i(\bsxi)$ are $n$-dimensional polynomials that are mutually orthogonal
with respect to the probability density function $\rho(\bsxi)$
\begin{equation}
\left<\Phi_i(\bsxi),\Phi_j(\bsxi)\right> = \int_\Xi \Phi_i(\bsxi) \Phi_j(\bsxi) \rho(\bsxi)\text{d}\bsxi =
||\Phi_i||^2 \delta_{ij},
\label{equ:ortrelation}
\end{equation}
with $\delta_{ij}$ the Kronecker delta
and $||\Phi_i||^2=\left<\Phi_i,\Phi_i\right>$. For independent random variables the $\Phi_i$'s are
tensor products of one-dimensional polynomials, $\Phi_i(\bsxi)=\prod_{j=1}^n \tilde\Phi_{j_i}(\xi_j)$.

To determine the polynomial chaos expansion coefficients $\bsxi_i$ and $\bsu_i(\bsx)$
there are two main-stream methods, viz. {\em Spectral Projection} and {\em Galerkin}. Both project onto 
the polynomial space, but
whereas in the spectral projection approach the {\em expansion} (\ref{equ:expansion}) is projected,
in the Galerkin approach the {\em governing equation} (\ref{equ:esto}), with the expansion substituted,  
is projected onto the polynomial space. Typically, for the random input variables the first approach
is chosen, resulting in
\begin{equation}
 \bsxi_j = \frac{\left<\bsxi,\Phi_j(\bsxi)\right>}{\left<\Phi_j(\bsxi),\Phi_j(\bsxi)\right>}=
 \frac{1}{||\Phi_j||^2}
 \int_\Xi \bsxi \, \Phi_j(\bsxi) \, \rho(\bsxi)\, \text{d}\bsxi,\qquad j=0,1, ... \;.
\label{eq:SPtheta}
\end{equation}
%since in most cases $\bsxi$ is a linear or almost linear function. 
If the function is truely nonlinear aliasing is a threat for the accuracy \cite{Xiu2007}, 
%Since, in this paper we want to concentrate on the
%essentials of our new method, i.e., the derivation of the polynomials for correlated multivariate distributions,
therefore we will use the Galerkin method to obtain the solution.
The expansions (\ref{equ:expansion})
are substituted into the stochastic equation (\ref{equ:esto}) before being projected
\begin{eqnarray}
&&\left<\mathcal{L}\left(\bsx,\sum_{i=0}^\infty \bsxi_i \Phi_i(\bsxi); \sum_{i=0}^\infty \bsu_i(\bsx)\Phi_i(\bsxi)\right),\Phi_j(\bsxi)\right>=\nonumber\\
&&\qquad\qquad\qquad\quad=\left<\bff\left(\bsx,\sum_{i=0}^\infty \bsxi_i \Phi_i(\bsxi)\right),\Phi_j(\bsxi)\right>,\qquad j=0,1,... \;.
\label{equ:Gu}
\end{eqnarray}
In practice, the number of expansion terms is truncated to 
\begin{equation}
 (N+1)=\frac{(n+p)!}{n!\,p!},
 \label{equ:dim}
\end{equation}
where $p$ is the highest order of the polynomials and the integrals are either computed exactly or
approximated using quadrature rules or Monte Carlo sampling. E.g., for a system of ODEs the differential operator in the left-hand side of (\ref{equ:esto})
reduces to $\dot{\bsu}(t,\bsxi)$ and its projection onto $\Phi_j$ to $\dot{\bsu}_j(t)\;||\Phi_j||^2$.

\subsection{Multivariate arbitrarily distributed input variables}
In \cite{G-S} a method has been proposed to extend the polynomial chaos method to an arbitrarily distributed 
univariate distribution. The one-dimensional polynomials $\Phi_i(\xi)$ are constructed to be mutual 
orthogonal with respect to an arbitrarily pdf $\rho(\xi)$ using the well-known and robust Gram-Schmidt 
orthogonalization method (see e.g. \cite{GL3}). This method can be extended to
multivariate {\em independent} random variables, where the orthogonal multidimensional polynomials are
the product of the constructed one-dimensional orthogonal polynomials. 
Here we use the same orthogonalization approach
to construct directly a multidimensional orthogonal polynomial basis, $\{\Phi_{j}(\bsxi)\}_{j=0}^{N}$,
for {\em correlated} multivariate random input variables.
Since the Gram-Schmidt method constructs from an arbitrary basis a basis that is orthogonal with respect to
a given innerproduct, the first step is the choice of a suitable set of linearly independent polynomials.
It should be noticed that any linearly independent set of polynomials can be used in this method,
but for simplicity we will use the set of monic polynomials $\{e_j(\bsxi)\}_{j=0}^{N}$
given by
\begin{equation}
e_j(\bsxi)=\prod_{l=1}^n\xi_l^{j_l},
\quad j=0, \ldots,N,\;j_l \in \{0,\ldots,p\},\text{ and }\sum_{l=1}^n j_l \leq p,
\label{equ:pol}
\end{equation}
%\JGB{[Can we make the sum more precise? I.e. incl. \# mon. of a certain degree]}
where $n$ is the dimension of the random vector $\bsxi$, $p$ the highest polynomial
degree chosen, and $N+1$ the dimension of the polynomial basis given by Equation (\ref{equ:dim}). 
E.g., if $n=p=2$ then $N+1=6$ and the set of linearly independent polynomials 
$\{e_{j}(\xi_{1},\xi_{2})\}_{j=0}^5$ equals $\{1, \xi_{1},\xi_{2},\xi_{1}^{2},\xi_{2}^{2}, \xi_{1}\xi_{2} \}$.
Next, the orthogonal polynomial basis $\{\Phi_{j}(\mbox{\boldmath$\xi$})\}_{j=0}^{N}$
is constructed  from $\{e_{j}(\mbox{\boldmath$\xi$})\}_{j=0}^{N}$ using the Gram Schmidt algorithm
\begin{eqnarray}
\Phi_{0}(\bsxi)&=&1, \nonumber\\
\Phi_{j}(\bsxi)&=&e_{j}(\bsxi)-\sum_{k=0}^{j-1}c_{jk}\Phi_{k}(\bsxi)\quad\text{for }1\leq j \leq N,
\label{equ:recur-formu}
\end{eqnarray}
where the coefficients $c_{jk}$ are given by
\begin{equation}
c_{jk}=\frac{\left<e_{j}(\mbox{\boldmath$\xi$}),\Phi_{k}(\mbox{\boldmath$\xi$})\right>}{\left<\Phi_{k}(\mbox{\boldmath$\xi$}),\Phi_{k}(\mbox{\boldmath$\xi$})\right>}
\label{equ:c_coeff}
\end{equation}
and the innerproduct is taken with respect to the pdf $\rho(\bsxi)$. Note, that the basis is not
unique, it is dependent on the choice and on the ordering of the set of polynomials $\{e_j\}$.

Since each polynomial of the orthogonal basis 
$\{\Phi_{j}(\mbox{\boldmath$\xi$})\}_{j=0}^{N}$
can be written as a sum of the monic 
polynomials $e_j(\bsxi)$, the inner products in Equation (\ref{equ:c_coeff}) can be calculated as sum of 
raw moments 
\begin{equation}
\mu^{r_{1}\cdots r_{n}}_{\bsxi}=\int_\Xi\prod_{l=1}^n\xi_{l}^{r_l}\rho(\bsxi)d\bsxi,\quad
r_l \in \{0,\ldots,2p\},\text{ and }\sum_{l=1}^n r_l \leq 2p.
\end{equation}
With this procedure we immediately get a set of multidimensional orthogonal polynomials. Note, that if the 
random variables are independent these polynomials are the same as the ones obtained by taking the tensor-products
of the one-dimensional orthogonal polynomials.
Once these multidimensional orthogonal polynomials have been computed the PCE method can be applied analogously
to the independent case, i.e., projection onto the polynomial space of either the expansion or 
the equation with the expansion substituted.

\subsubsection{Statistics}
Part of the ease of use of PCE is the simplicity with which one obtains the most used statistics of the
QoI: mean, variance, and Sobol' indices can be directly expressed using the expansion coefficients \cite{Sudret2008}.
As a consequence of the orthogonality of the basis $\{\Phi_{j}(\mbox{\boldmath$\xi$})\}_{j=0}^{N}$
this favorable feature still holds 
for the mean and variance , i.e., mean and variance of the solution vector $\bsu$ in 
Equation (\ref{equ:expansion}) are given by
\begin{eqnarray}
 \bsmu_\bsxi(\bsx) &=& \int_\Xi \bsu(\bsx,\bsxi) \rho(\bsxi) d\bsxi = \bsu_0(\bsx),\text{ and}\\
 \bssi^2_\bsxi(\bsx) &=& \int_\Xi \left(\bsu(\bsx,\bsxi)-\bsmu_\bsxi(\bsx)\right)^2 \rho(\bsxi) d\bsxi =
 \sum_{i=1}^\infty \bsu^2_i(\bsx)||\Phi_i(\bsxi) ||^2
 .
\end{eqnarray}
The Sobol' indices \cite{Sobol1990,Li-globalsa_2010, Sudret2008, Owen2013} measure the influence of varying
{\em only the specified combination of} variables on the total variance. 
They are based on the terms of the Sobol' decomposition of the QoI
\begin{equation}
 \bsu(\bsx,\bsxi) = \bsM_0(\bsx) + \sum_{\bsl \subseteq \{1,...,n\}} \bsM_\bsl(\bsx,\bsxi_\bsl),
 \label{eq:Soboldec}
\end{equation}
where $\bsM_0(\bsx) = \bsmu_\bsxi(\bsx)$, and the terms
$\bsM_\bsl$ are recursively given by
\begin{equation}
 \bsM_\bsl(\bsx,\bsxi_\bsl) = \E_{\bsxi_{- \bsl}}\left[\bsu|\bsxi_\bsl\right] - 
 \left(\bsM_0(\bsx) + \sum_{\bsk \subset \bsl} \bsM_\bsk(\bsx,\bsxi_\bsk)\right),
\end{equation}
with  the marginal expectation given by
\begin{equation}
\E_{\bsxi_{-\bsl}}\left[\bsu(\bsx,\bsxi)|\bsxi_{\bsl}\right] =
 \int_{\Xi_{-\bsl}} \bsu(\bsx,\bsxi) 
 \left(\int_{\Xi_\bsl} \rho(\bsxi) d\bsxi_{\bsl}\right)
 d\bsxi_{-\bsl},
\end{equation}
where $\bsxi_{-\bsl}$ indicates all elements of $\bsxi$ except $\bsxi_\bsl$. It is easy to see that this
expansion holds, since $\bsM_{1,...,n}(\bsx,\bsxi)$ equals $\bsu(\bsx,\bsxi)$ minus all previous terms.
Note, that this is an extension to general distributions of the original definition of Sobol' \cite{Sobol1990},
which was for $\mathcal{U}[0,1]^n$.

Using decomposition (\ref{eq:Soboldec}) the variance of the QoI is decomposed
\begin{eqnarray}\label{eq:vardec}
&& \Var\left[\bsu\right] =\Cov\left[\bsu,\bsu\right]= \Cov[\bsM_0 + \sum_{\bsl \subseteq \{1,...,n\}} \bsM_\bsl,\bsu] = \nonumber\\
&&  \sum_{\bsl \subseteq \{1,...,n\}} \Cov\left[\bsM_\bsl,\bsu\right] 
= \sum_{\bsl \subseteq \{1,...,n\}} \left(\Var\left[\bsM_\bsl \right] + \Cov\left[\bsM_\bsl,\bsu - \bsM_\bsl \right]\right).
\end{eqnarray}

The Sobol' indices are now defined as the - normalized - covariance between the respective terms in the expansion 
and the QoI 
\begin{equation}
 \bsS_\bsl(\bsx) = \frac{\Cov\left[\bsM_\bsl(\bsx,\bsxi_\bsl),\bsu\right]}
 {\Var\left[\bsu\right]},\quad\text{for }\bsl \subseteq \{1,...,n\}.
\end{equation}

Following \cite{Li-globalsa_2010} (cf. also Eq. (\ref{eq:vardec})), 
we decompose the Sobol' indices in a part that is dependent on the respective
variables and a part that is due to the correlation with all other variables
\begin{eqnarray}
 \bsS_\bsl(\bsx) &=&  \bsS_\bsl^u(\bsx) + \bsS_\bsl^c(\bsx), \text{ with}\label{eq:Si}\\
 \bsS_\bsl^u(\bsx) &=& \frac{\Var\left[\bsM_\bsl(\bsx,\bsxi_\bsl)\right]}
 {\Var\left[\bsu\right]}, \text{ and}\label{eq:Si_u}\\
 \bsS_\bsl^c(\bsx) &=& \frac{\Cov\left[\bsM_\bsl(\bsx,\bsxi_\bsl),\bsu-\bsM_\bsl(\bsx,\bsxi_\bsl)\right]}
 {\Var\left[\bsu\right]};\label{eq:Si_c}
\end{eqnarray}
where the index $\bsS_\bsl(\bsx)$ represents that part of
the variance of the QoI that is due to the variance of the set of variables $\bsxi_{\bsl}$, 
$\bsS_\bsl^u(\bsx)$ represents the uncorrelated 
share of it,
that is the contribution to the variance that comes from the set of variables $\bsxi_{\bsl}$ 
by themselves, and $\bsS_\bsl^c(\bsx)$ represents the correlated share, 
the contribution to the variance of the QoI that comes from 
the correlation of the set of the variables $\bsxi_{\bsl}$ with 
the set of variables $\bsxi_{-_\bsl}$.
The {\em total influence} of a single variable, $\xi_{j}$, including that part of the variance due to 
variable $\xi_{j}$ alone
and the fraction due to any combination of $\xi_{j}$ with the remaining variables is given by
\begin{equation}
 \bsS_{\bs T_{j}}(\bsx) = \sum_{\bsl\subseteq \{1,...,n\} \wedge j\in\bsl} \bsS_\bsl(\bsx),
\end{equation}
and analogous definitions for $\bsS_{\bs T_{j}}^{u}(\bsx)$, and $\bsS_{\bs T_{j}}^{c}(\bsx)$.

By definition $\bsS_\bsl^u(\bsx)$ is always positive, but $\bsS_\bsl^c(\bsx)$ can be either
positive or negative and therefore $\bsS_\bsl(\bsx)$ can also be either positive or negative. 
This makes the interpretation and more specific the influence-based ranking of the variables less easy 
than in the case of independent variables. For the variables with small $\bsS^c$ values the ``normal ranking''
can still be used, but when the Sobol' indices for a variable have a large $\bsS^c$ part, the interpretation 
is still an open issue.

\medskip 
For the computation of the Sobol' indices it is useful to notice that
a Sobol' index can be written as a 
sum of Sobol' indices of the $\Phi$-terms with the appropriate coefficient,
since $\bsu(\bsx,\bsxi)$ is a linear combination of $\Phi$'s (see also \cite{Sudret2008}).
For independent random variables the polynomials $\Phi$ (and the distribution $\rho$) are products of
univariate contributions which implies that 

\begin{eqnarray}
 \E_{\bsxi_{- \bsl}}\left[\Phi_j(\bsxi) | \bsxi_\bsl\right] &=& 
 \int_{\Xi_{- \bsl}}  \prod_{k=1}^n \tilde\Phi_{j_k}(\xi_k)
 \left(\int_{\Xi_\bsl} \rho(\bsxi) d\bsxi_\bsl\right)
 d\bsxi_{- \bsl} =\nonumber\\
 &=& \prod_{k \in \bsl} \tilde\Phi_{j_k}(\xi_k) 
 \int_\Xi \prod_{k \in -\bsl} \tilde\Phi_{j_k}(\xi_k) \rho(\bsxi)d\bsxi =\nonumber\\
 &=& \left\{ \begin{array}{ll}
              \prod_{k \in \bsl} \tilde\Phi_{j_k}(\xi_k) & 
              \text{if }\prod_{k \in -\bsl} \tilde\Phi_{j_k}(\xi_k)=1(=\tilde\Phi_0),\\
              0 & \text{otherwise.}
             \end{array}
      \right.
 \label{eq:CondE_indep}
\end{eqnarray}
For a non-zero marginal expectation, the contribution to the numerator of the Sobol' index, consisting only of the 
uncorrelated part (\ref{eq:Si_u}), is then given by
\begin{equation}
 \Var\left[\E_{\bsxi_{-_\bsl}}\left[\Phi_j(\bsxi)|\bsxi_\bsl\right]\right] =
 \int_\Xi  \left(\prod_{k \in \bsl} \tilde\Phi_{j_k}(\xi_k)\right)^2 \rho(\bsxi)d\bsxi
 =||\prod_{k \in \bsl} \tilde\Phi_{j_k}||^2.
\end{equation}
So, if the polynomial
\begin{equation*}
 \Phi_j(\bsxi) = \prod_{k \in \bsl} \tilde\Phi_{j_k}(\xi_k) \prod_{k \in -\bsl} \tilde\Phi_0(\xi_k),
\end{equation*}
the contribution to the Sobol' index is given by the expansion coefficient $u_j(\bsx)$ 
and otherwise it is zero.

For correlated random variables the orthogonality property used in (\ref{eq:CondE_indep}) no longer holds
and a Sobol' index is no longer given by a simple combination of expansion coefficients. 
Fortunately, the same trick of separating the variables can still be used to compute the marginal 
expectation and the covariance thereof, albeit at a lower level: 
each $\bsu(\bsx,\bsxi)$ is a linear combination of $\Phi$'s which themselves are linear combinations
of the monic polynomials $e_j(\bsxi)$, so each Sobol' index can be calculated as a linear combination of 
already computed raw moments with respect to the full pdf
\begin{eqnarray}
 \E_{\bsxi_{- \bsl}}\left[e_j(\bsxi) | \bsxi_\bsl\right] &=& 
 \int_{\Xi_{- \bsl}}  \prod_{k=1}^n \xi_k^{j_k}
 \left(\int_{\Xi_\bsl} \rho(\bsxi) d\bsxi_\bsl\right)
 d\bsxi_{- \bsl}\nonumber\\
 &=& \prod_{k \in \bsl} \xi_k^{j_k} 
 \int_\Xi \prod_{k \in -\bsl} \xi_k^{j_k} \rho(\bsxi)d\bsxi\nonumber\\
 &=& \mu^{j_{-\bsl}}_{\bsxi} \prod_{k \in \bsl} \xi_k^{j_k}.
 \label{eq:CondE_corr}
\end{eqnarray}
Note, that in this short-hand notation 
$\mu^{i_{\bsl}}_{\bsxi}\;\mu^{j_{\bsk}}_{\bsxi} \neq \mu^{i_{\bsl}+j_{\bsk}}_{\bsxi}$.\\
The $e_j$-contribution to the numerator of the Sobol index $\bsS_\bsl^u$ (\ref{eq:Si_u}) is then given by
\begin{eqnarray}
 \Var\left[\E_{\bsxi_{-_\bsl}}\left[e_j(\bsxi)|\bsxi_\bsl\right]\right] &=&
 (\mu^{j_{-\bsl}}_{\bsxi})^2 \left(
 \int_\Xi \prod_{k \in \bsl} \xi_k^{2 j_k} \rho(\bsxi)d\bsxi
 -\left(\int_\Xi \prod_{k \in \bsl} \xi_k^{j_k} \rho(\bsxi)d\bsxi\right)^2
 \right)\nonumber\\
 &=& (\mu^{j_{-\bsl}}_{\bsxi})^2 [\mu^{2 j_{\bsl}}_{\bsxi}-(\mu^{j_{\bsl}}_{\bsxi})^2].
\end{eqnarray}
Analogously, the numerator of the Sobol' index $\bsS_\bsl^c$ (\ref{eq:Si_c}) can be written as a combination of
moments, since
\begin{eqnarray}
 \lefteqn{\Cov\left[\E_{\bsxi_{-_\bsl}}\left[e_j(\bsxi)|\bsxi_\bsl\right],e_i(\bsxi)\right]}\nonumber\\
 \qquad&=&
 \mu^{j_{-\bsl}}_{\bsxi}
 \int_{\Xi} \prod_{k \in \bsl} \xi_k^{j_k} \prod_{k=1}^n \xi_k^{i_k}
 \rho(\bsxi) d\bsxi - (\mu^{j_{-\bsl}}_{\bsxi} \mu^{j_{\bsl}}_{\bsxi}) \mu^{i_{1...n}}_{\bsxi}
 \nonumber\\
 &=& \mu^{j_{-\bsl}}_{\bsxi} \mu^{j_{\bsl}+i_{1...n}}_{\bsxi}-(\mu^{j_{-\bsl}}_{\bsxi} \mu^{j_{\bsl}}_{\bsxi}) \mu^{i_{1...n}}_{\bsxi}.
\end{eqnarray}
All these moments have already been computed in the construction of the polynomial basis 
$\{\Phi_{j}(\mbox{\boldmath$\xi$})\}_{j=0}^{N}$, which
implies that also the Sobol' indices can be obtained without additional cost, although some bookkeeping is
required.

\bigskip
{\bf Remark.} Note, that the method is applicable to any type of continuous or discrete input probability distribution including
an experimentally obtained one. The only requirement is that the innerproducts of monic polynomials can 
be calculated, which are moreover only dependent on the distribution and not on the problem at hand.

\section{Examples}
We illustrate the proposed method with two examples to demonstrate the
significant influence of correlations between stochastic inputs on the distribution of the QoI.
The first example is a scalar ODE in two-dimensional random space; for this example we study also its
numerical behavior.
The second example is a set of four ODEs in three-dimensional random space, where 
we show the propagation of the resulting uncertainty in the parameters into the uncertainty in the QoI,
the concentration of the product.

\subsection{Decay equation}
Consider the scalar ODE
\begin{equation}
y'(t; \alpha,\beta) = -\alpha(y(t; \alpha,\beta)-\beta), \qquad 0 \leq t\leq 1, \qquad y(0) = 0,
\label{equ:Example1}
\end{equation}
with two random jointly distributed input variables $(\alpha, \beta) \sim \mathcal{N}(\bsmu,\bsSi)$, 
where the mean $\bsmu$ and the covariance matrix $\bsSi$ are given by
\begin{equation}
 \bsmu=\left(\begin{array}{c}1\\1\end{array}\right) \text{ and }
 \bsSi=\left(\begin{array}{c c} \sigma_\alpha^2 & \varrho\sigma_\alpha\sigma_\beta\\
                                \varrho\sigma_\alpha\sigma_\beta & \sigma_\beta^2
             \end{array}\right) \text{, with }
 \sigma_\alpha=\sigma_\beta=0.25.
\end{equation}
The correlation coefficient $\varrho$ varies from uncorrelated to fully correlated, i.e.,
$\varrho = \{0,\pm 0.5,\pm 0.9,\pm 1\}$.
Note, that the random input to the problem is univariate when the correlation coefficient $\varrho=\pm1$
\begin{equation}
y'(t; \alpha)=\left\{
\begin{array}{ll}
-\alpha(y(t; \alpha)-\alpha),&\varrho=+1\\
-\alpha(y(t; \alpha)+\alpha-2),&\varrho=-1
\end{array}
\right\}\quad
\begin{array}{l}
0 \leq t\leq 1,\; y(0)=0,\\
\alpha \sim \mathcal{N}(1,0.0625).
\end{array}
%\right.
\label{equ:Example1DPN}
\end{equation}

We want to show the effects of correlation and the convergence rate of the method for increasing expansion 
order uncontaminated with errors in the computation of the polynomials or the projections. 
Since the analytic expression of the  pdf of the correlated Gaussian distribution is known,
all moments, and the integrals thereof - necessary to compute the polynomial basis, the projection of the 
right hand-side, and the statistics - are computed exactly with the moment-generating function using the Symbolic toolbox of Matlab \cite{Matlab}.
The projection of the truncated Eq. (\ref{equ:Gu}) onto the polynomial basis results in a system of ODEs 
of the size of the expansion order. This system is solved with the accurate Matlab ODE solver ODE45 with
AbsTol 1e-6 (default) and RelTol changed to 1e-6. 

\begin{figure}[ht!]
   \centering
	 \subfloat[order 0, $\varrho=-0.5$]{
        \label{fig:Neg0}         
        \includegraphics[width=0.31\textwidth]{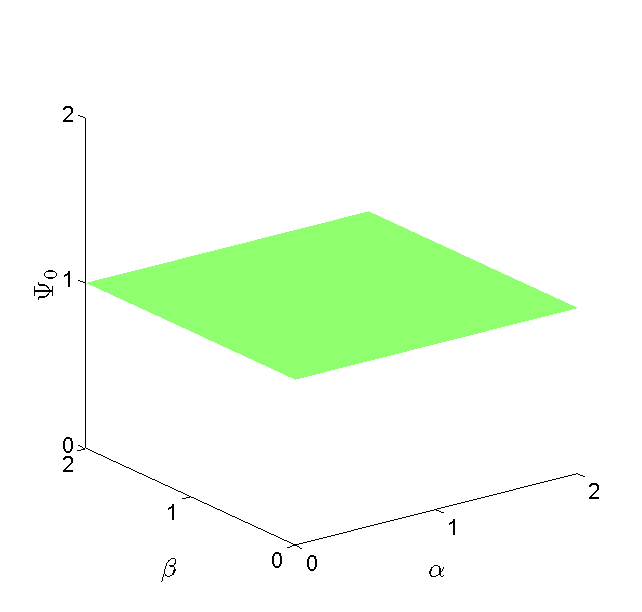}}
   \hspace{0.001\linewidth}
   %%----segunda subfigura----
   \subfloat[order 1, $\varrho=-0.5$]{
        \label{fig:Neg12}         
        \includegraphics[width=0.31\textwidth]{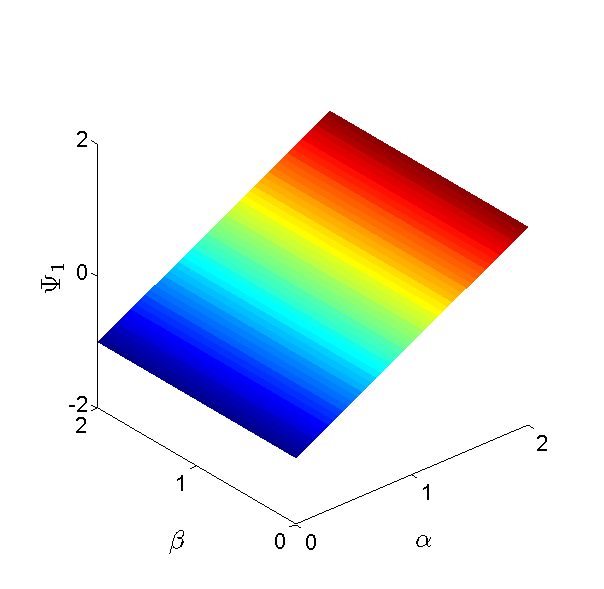}
         \includegraphics[width=0.31\textwidth]{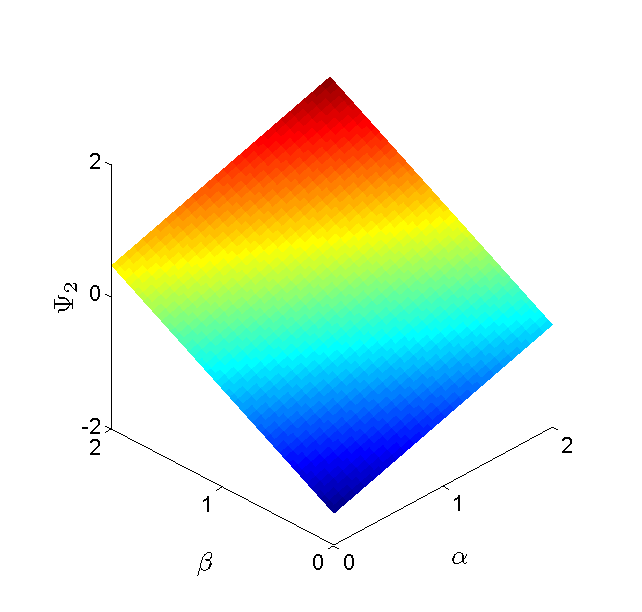}}\\
         
  	 \subfloat[order 0, $\varrho=0$]{
        \label{fig:U0}         
        \includegraphics[width=0.31\textwidth]{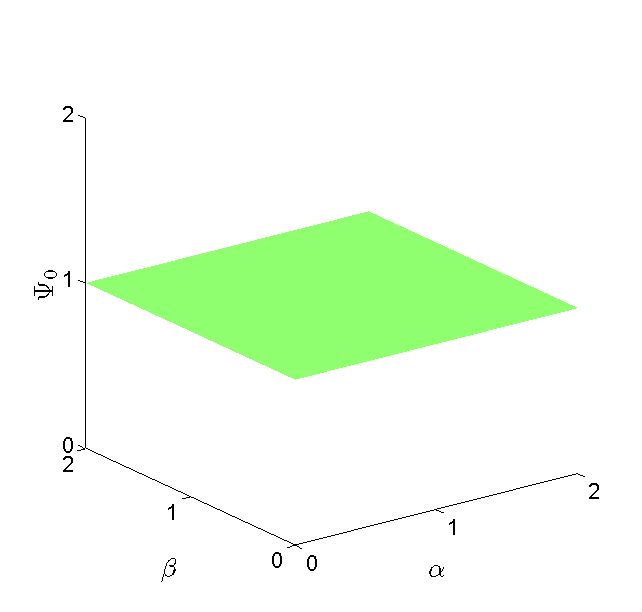}}
   \hspace{0.001\linewidth}
   %%----segunda subfigura----
   \subfloat[order 1, $\varrho=0$]{
        \label{fig:U12}         
        \includegraphics[width=0.31\textwidth]{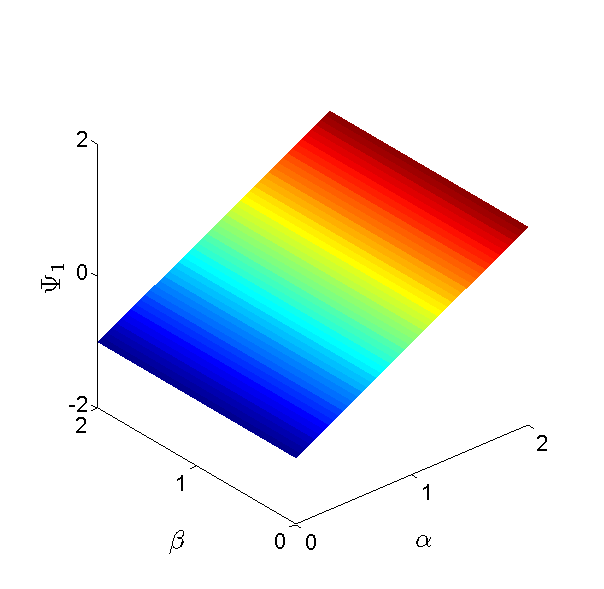}
         \includegraphics[width=0.31\textwidth]{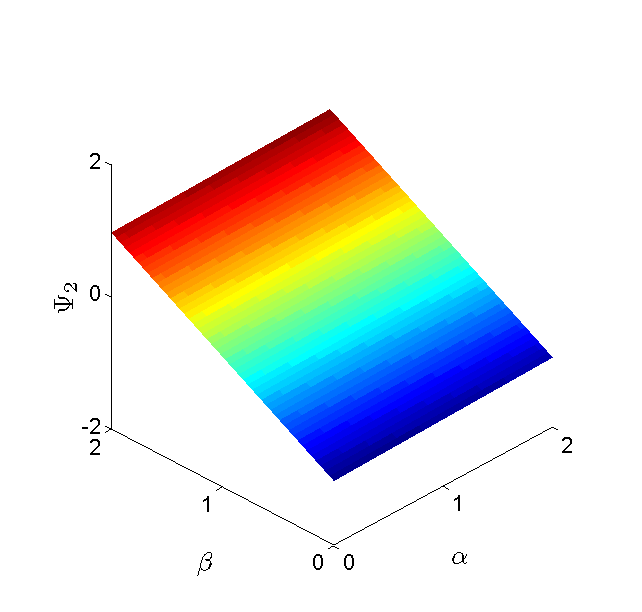}}\\

  \subfloat[order 0, $\varrho=0.5$]{
        \label{fig:Pos0}         
        \includegraphics[width=0.31\textwidth]{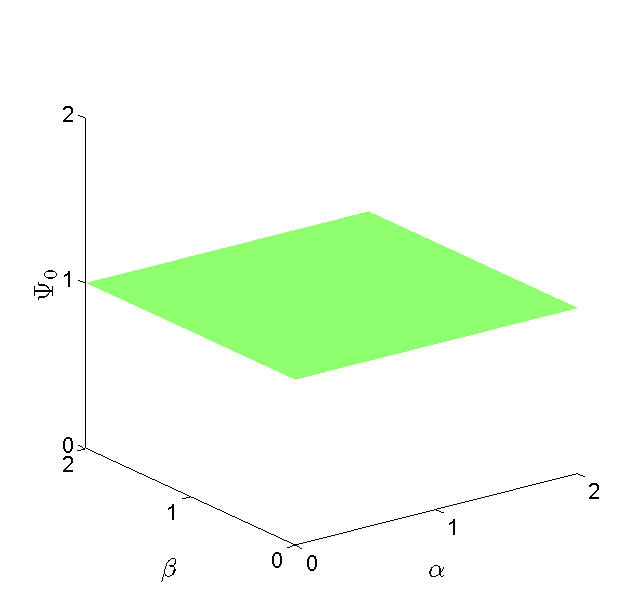}}
   \hspace{0.001\linewidth}
   %%----segunda subfigura----
   \subfloat[order 1, $\varrho=0.5$]{
        \label{fig:Pos12}         
        \includegraphics[width=0.31\textwidth]{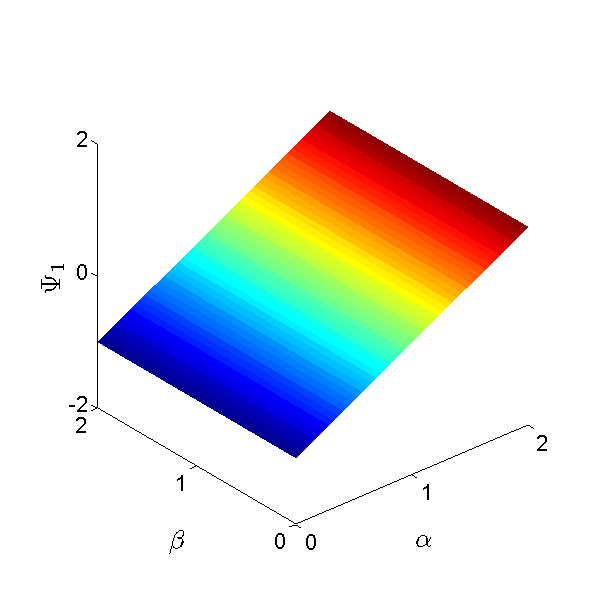}
         \includegraphics[width=0.31\textwidth]{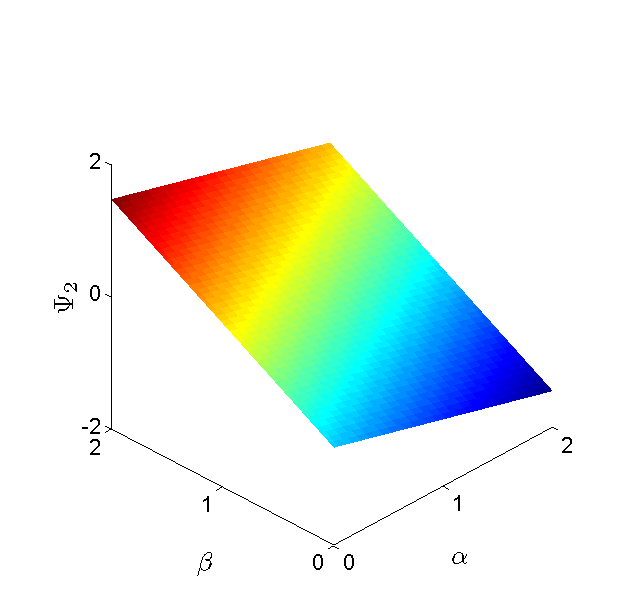}}\\

   \caption{Polynomial basis. Order 0 and 1, for $\varrho = \{0,\pm 0.5\}$}
  \label{fig:Pol}
\end{figure}

\begin{figure}[ht!]
   \centering
	 \subfloat{
        \label{fig:Neg23}         
        \includegraphics[width=0.31\textwidth]{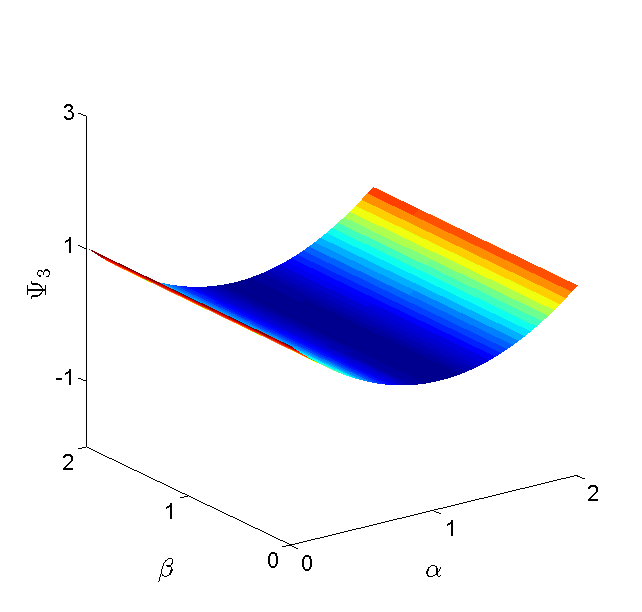}}
   \hspace{0.001\linewidth}
   %%----segunda subfigura----
        \subfloat[order 2, $\varrho=-0.5$]{
        \label{fig:Neg24}         
        \includegraphics[width=0.31\textwidth]{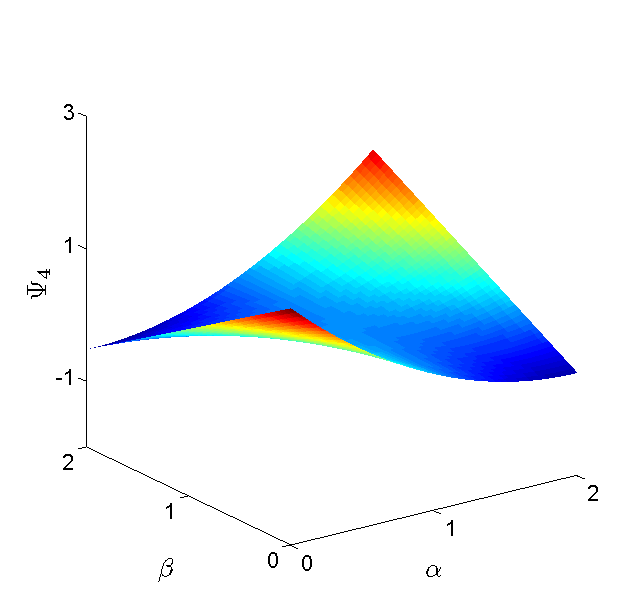}}
   \hspace{0.001\linewidth}
   \subfloat{
        \label{fig:Neg25}         
        \includegraphics[width=0.31\textwidth]{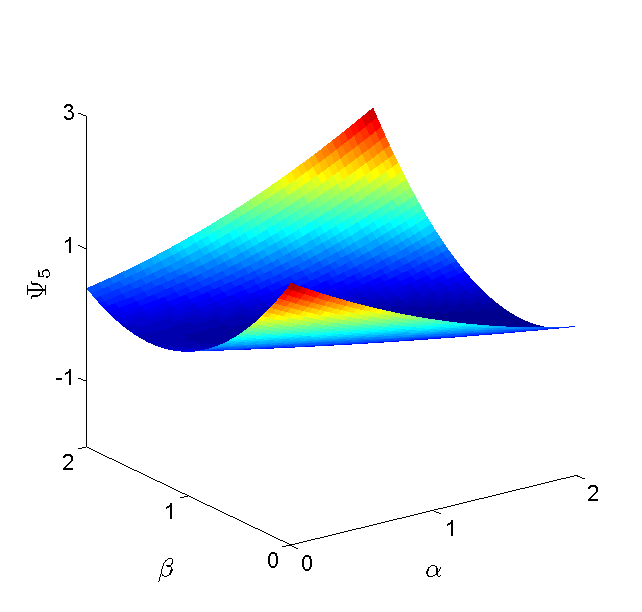}}\\
        	 \subfloat{
        \label{fig:U23}         
        \includegraphics[width=0.31\textwidth]{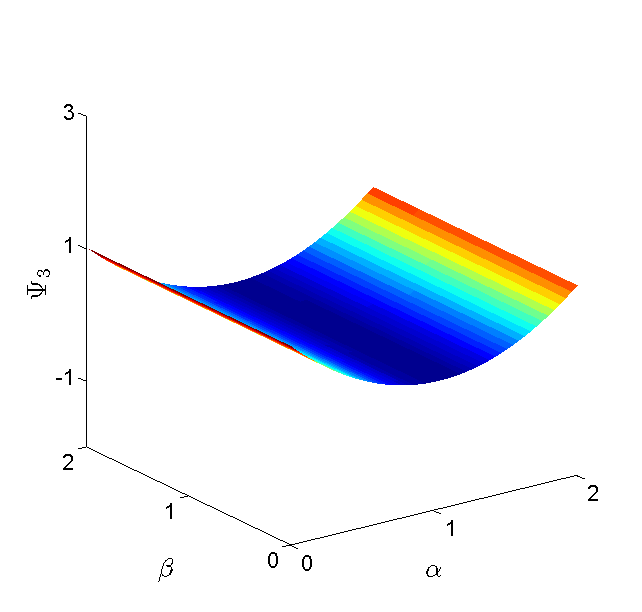}}
   \hspace{0.001\linewidth}
   %%----segunda subfigura----
        \subfloat[order 2, $\varrho=0$]{
        \label{fig:U24}         
        \includegraphics[width=0.31\textwidth]{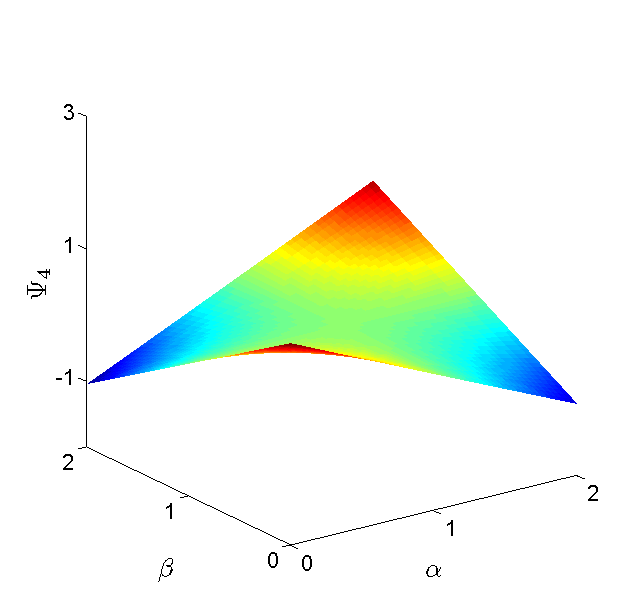}}
   \hspace{0.001\linewidth}
   \subfloat{
        \label{fig:U25}         
        \includegraphics[width=0.31\textwidth]{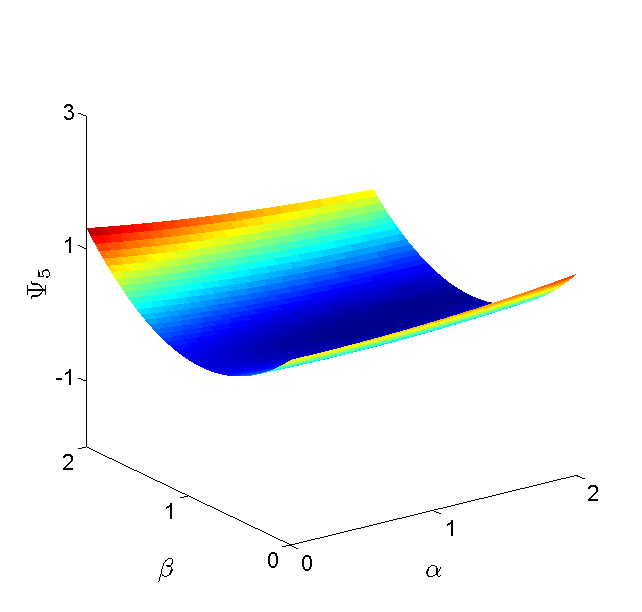}}\\
        	 \subfloat{
        \label{fig:Pos23}         
        \includegraphics[width=0.31\textwidth]{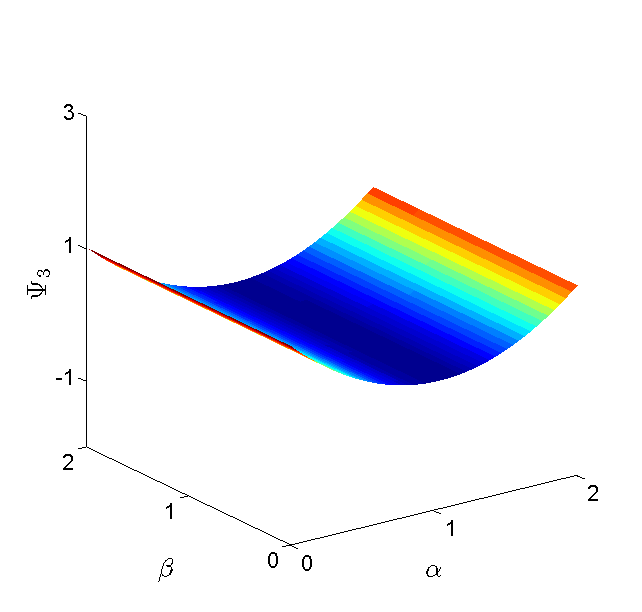}}
   \hspace{0.001\linewidth}
   %%----segunda subfigura----
        \subfloat[order 2, $\varrho=+0.5$]{
        \label{fig:Pos24}         
        \includegraphics[width=0.31\textwidth]{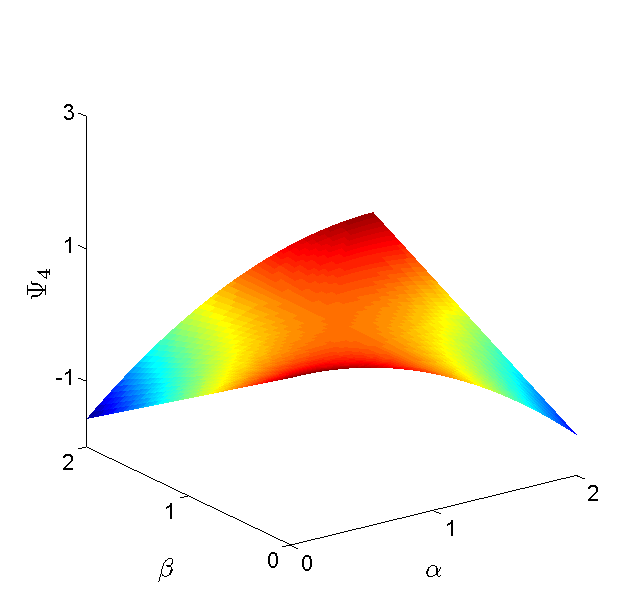}}
   \hspace{0.001\linewidth}
   \subfloat{
        \label{fig:Pos25}         
        \includegraphics[width=0.31\textwidth]{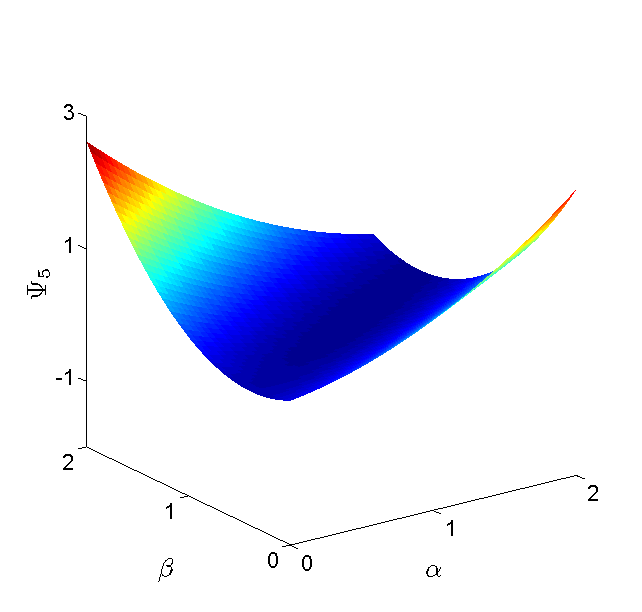}}\\
   \caption{Polynomial basis. Order 2, for $\varrho = \{0,\pm 0.5\}$}
  \label{fig:Pol2}
\end{figure}

\begin{figure}[htbp!]
   \centering
   \subfloat[Mean]{
        \includegraphics[width=6.5cm,height=6.5cm]{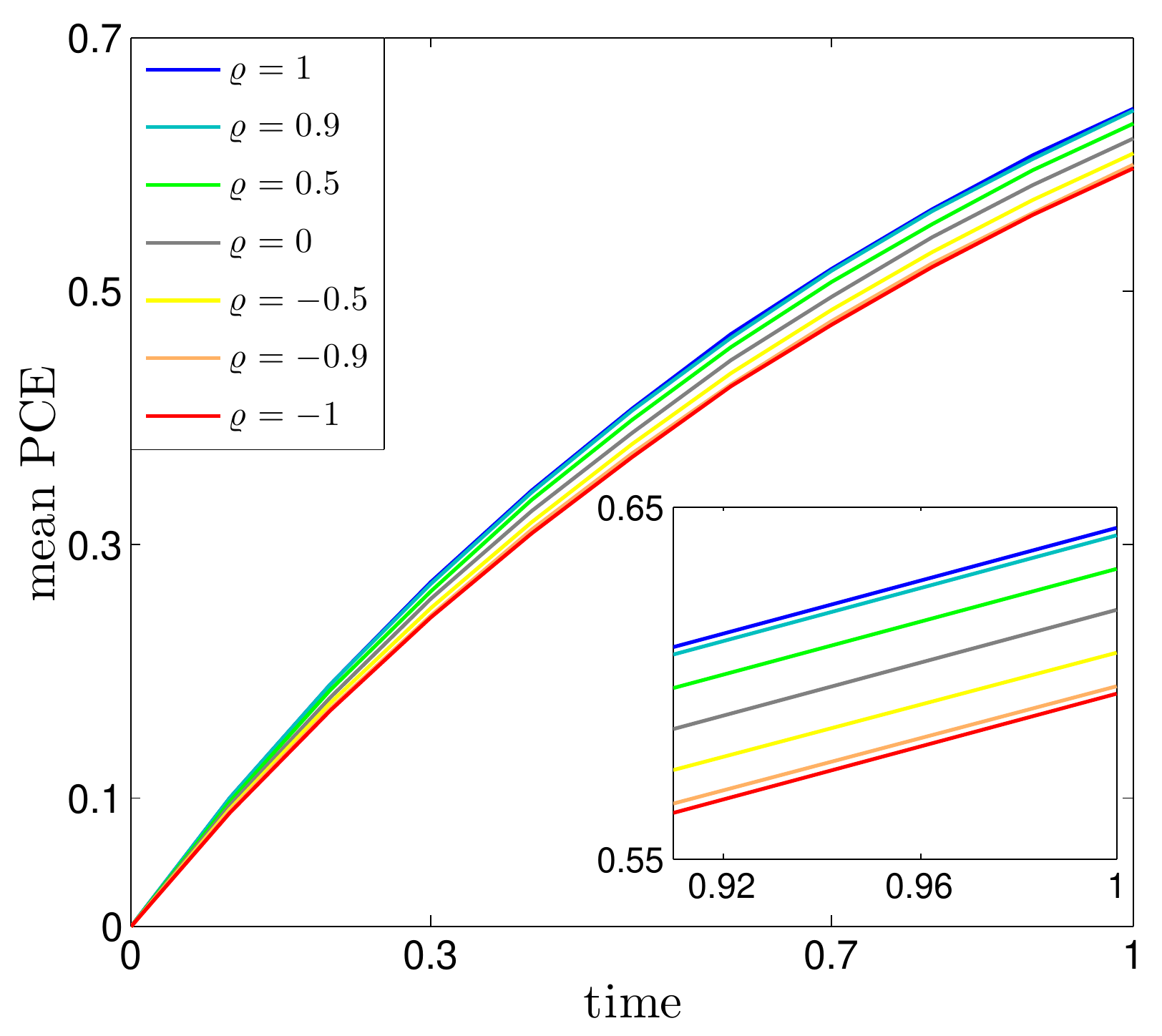}}
   \hspace{0.01\linewidth}
   \subfloat[Standard deviation]{
        \includegraphics[width=6.5cm,height=6.5cm]{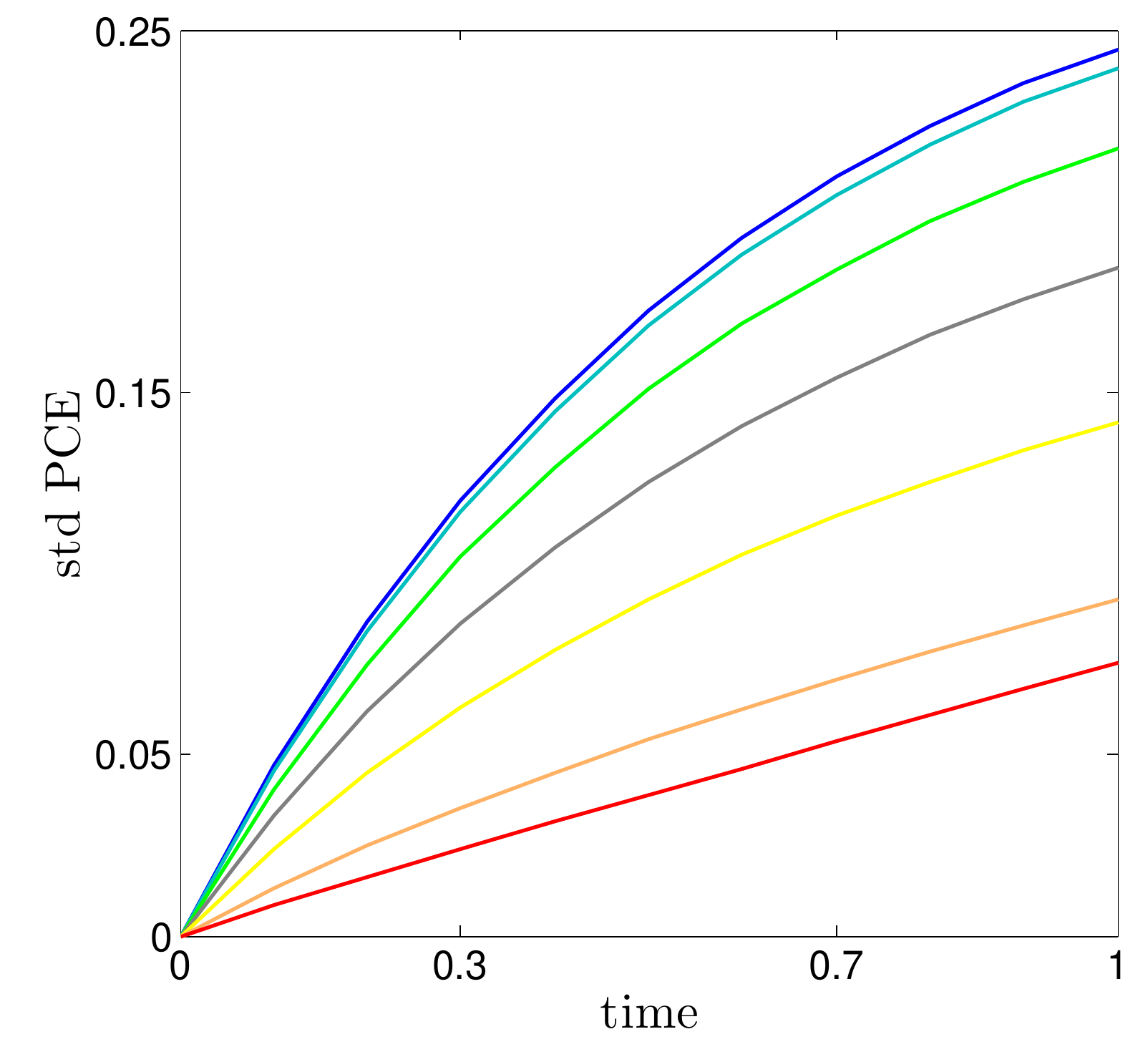}}
   \hspace{0.01\linewidth}
	\caption{Effect of correlation on mean and standard deviation of decay equation (\ref{equ:Example1})}
   \label{fig:Ex1-MSd} 
\end{figure}

\begin{figure}[htbp!]
\begin{tabular}{cccc}
\includegraphics[scale=0.21]{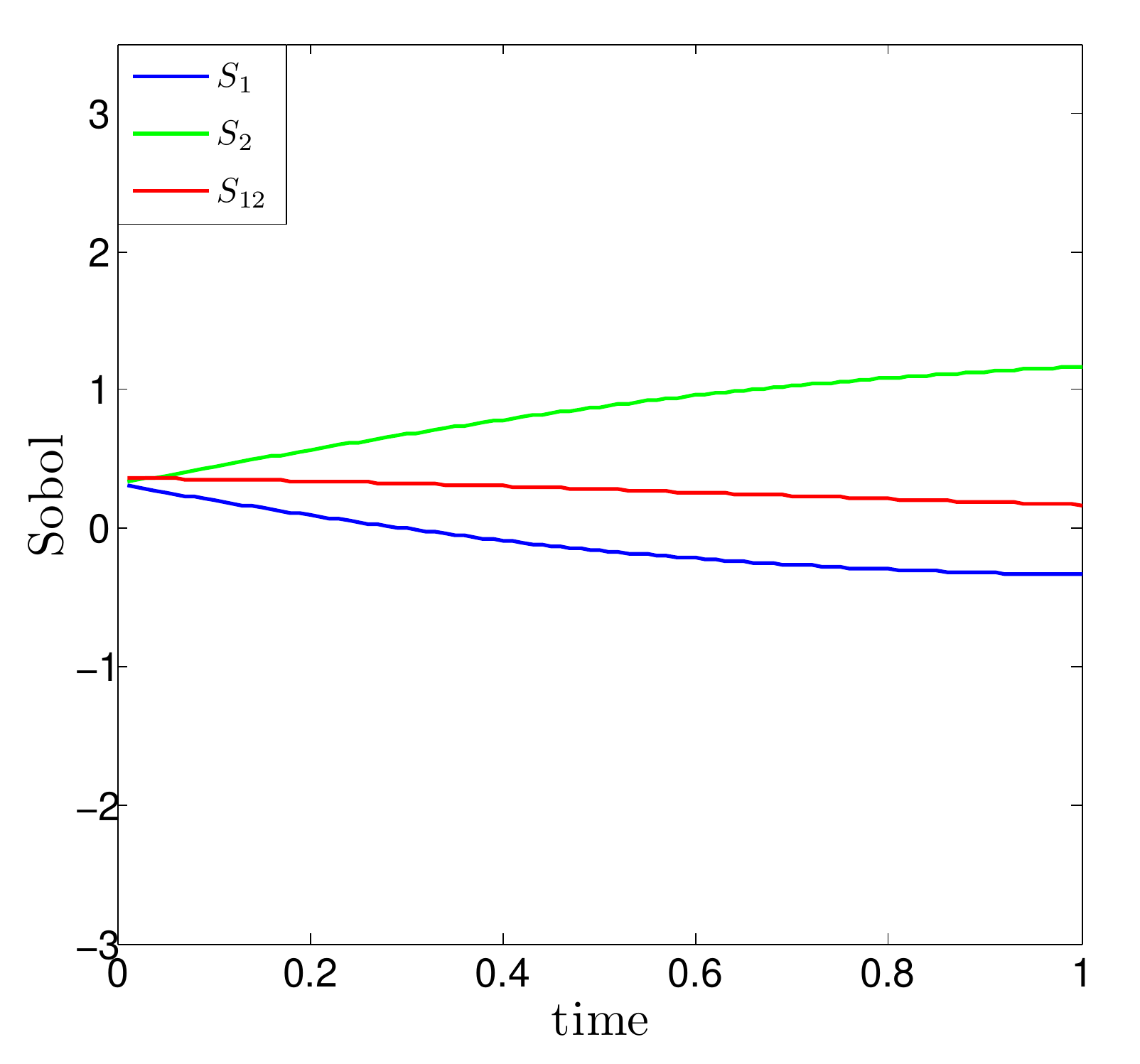}
\includegraphics[scale=0.21]{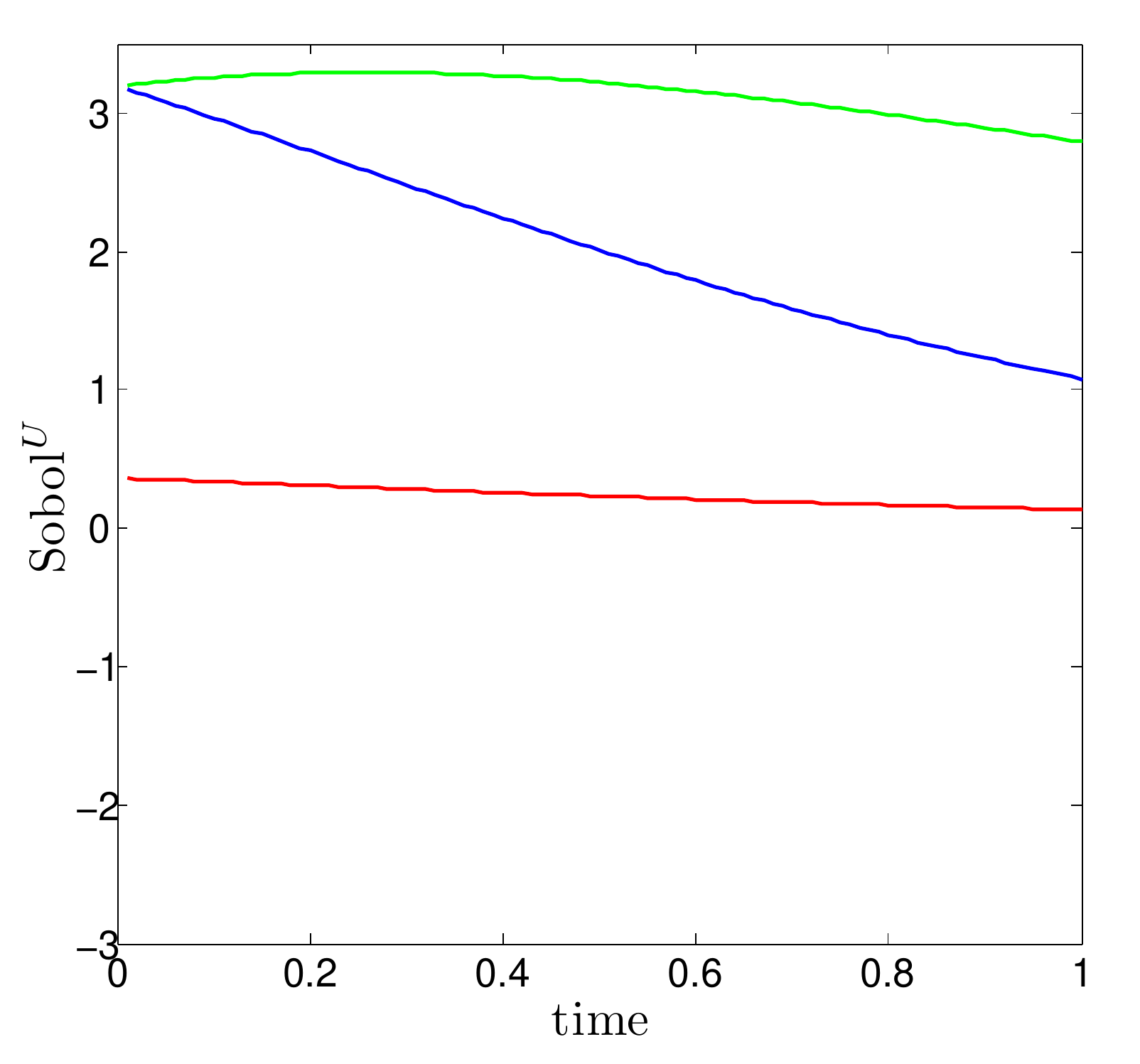}
\includegraphics[scale=0.21]{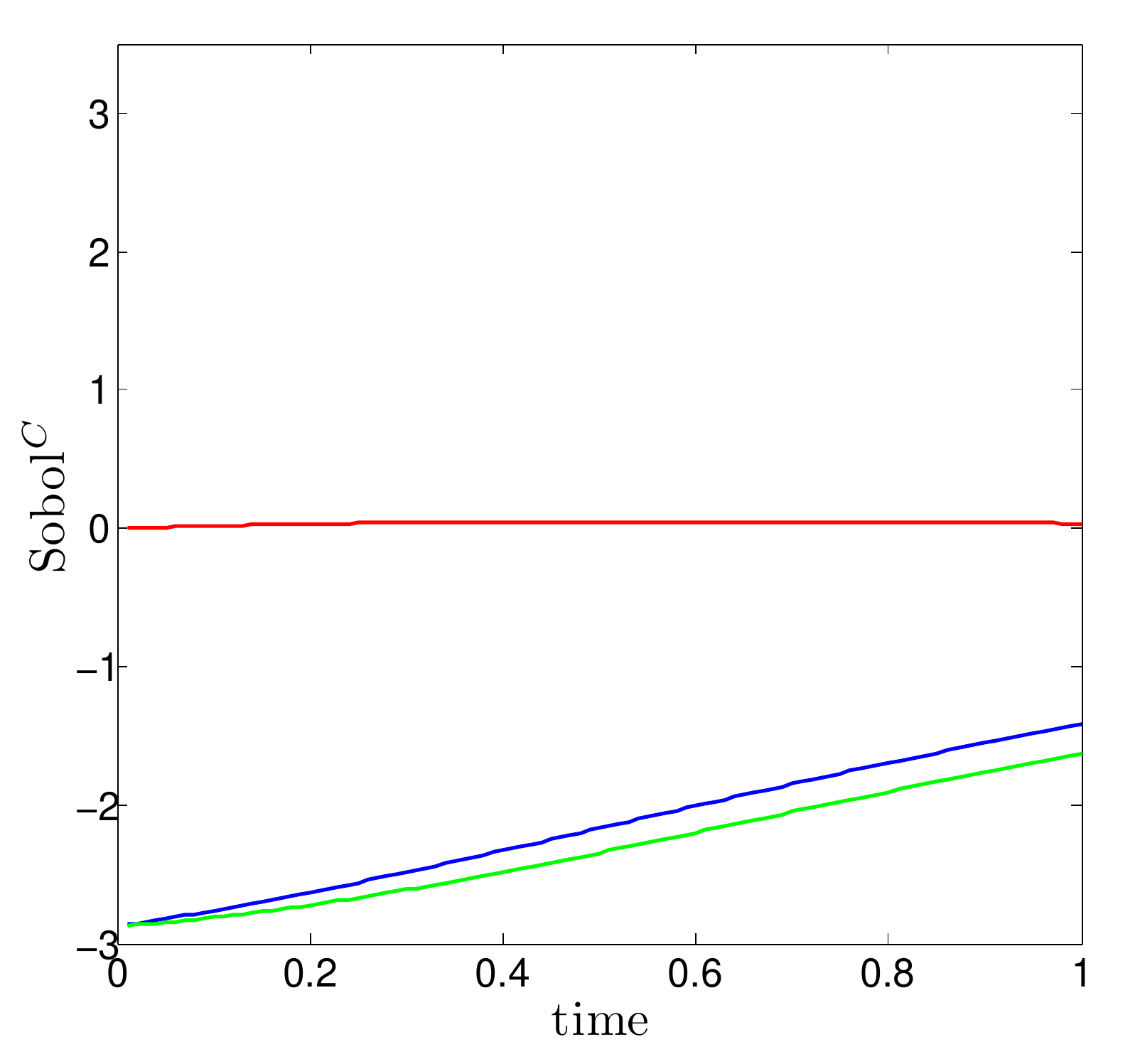}\parbox{2cm}{\vspace{-3cm} $\varrho=-0.9$}\\
\includegraphics[scale=0.21]{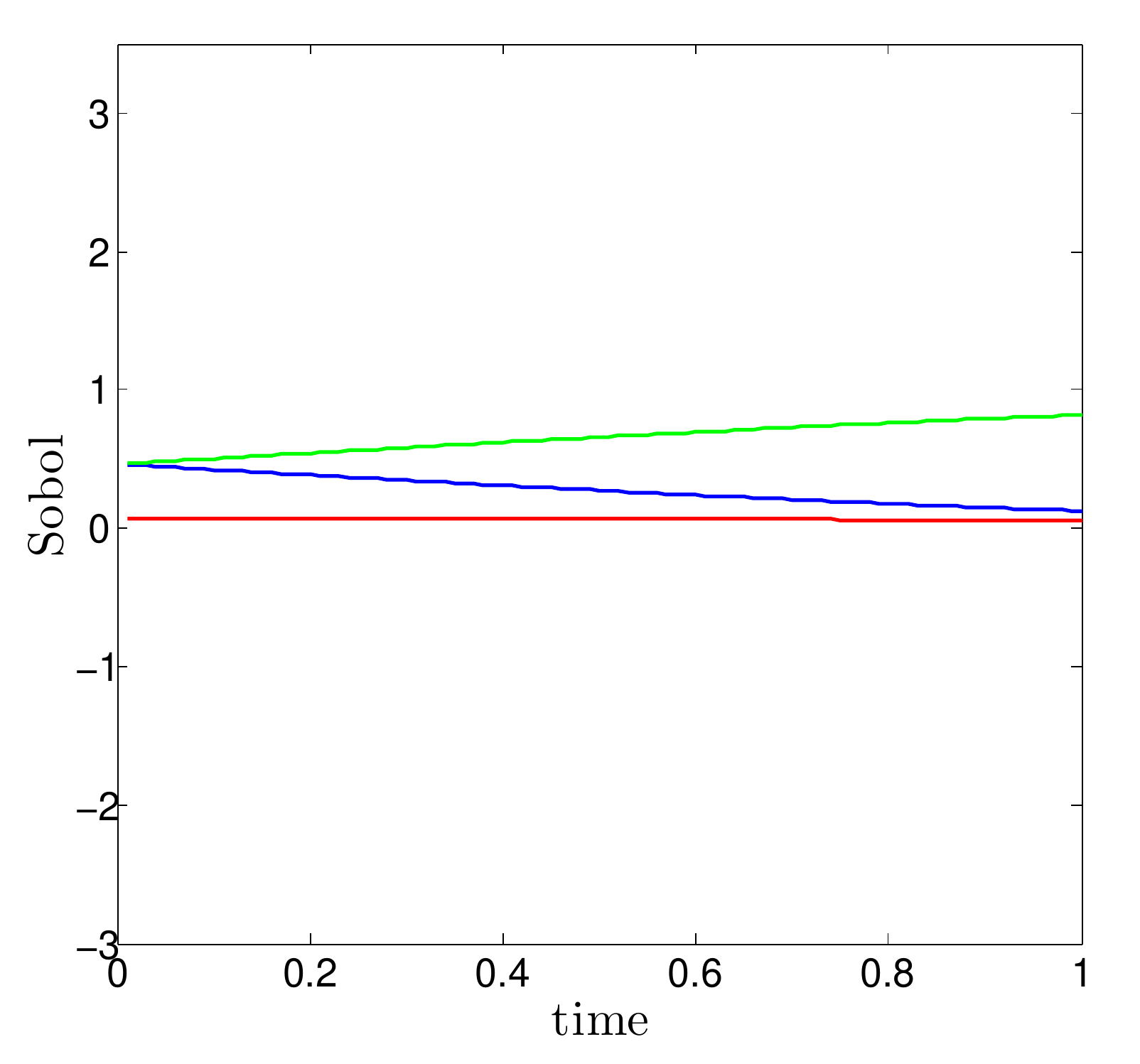}
\includegraphics[scale=0.21]{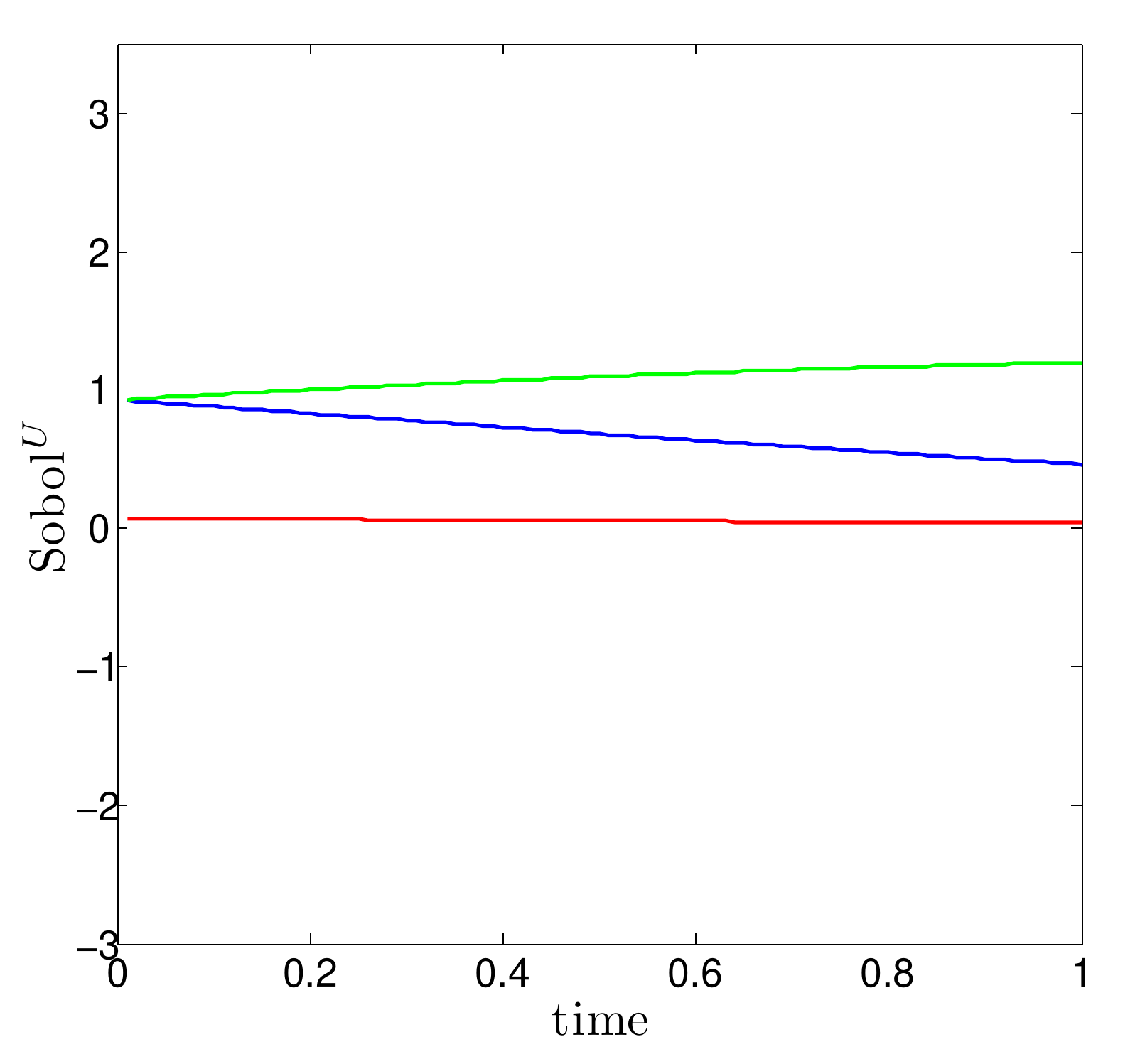}
\includegraphics[scale=0.21]{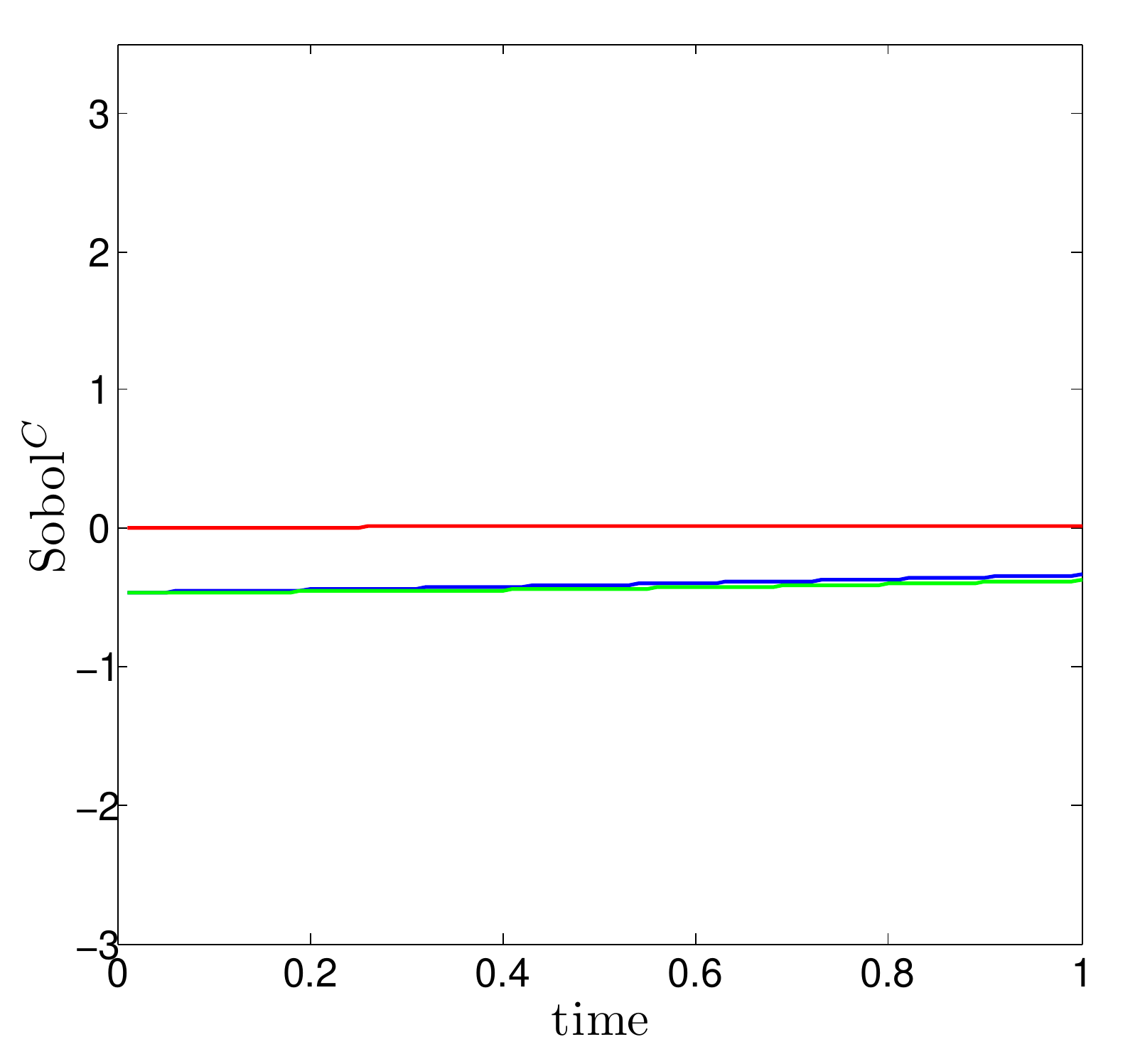}\parbox{2cm}{\vspace{-3cm} $\varrho=-0.5$}\\
\includegraphics[scale=0.21]{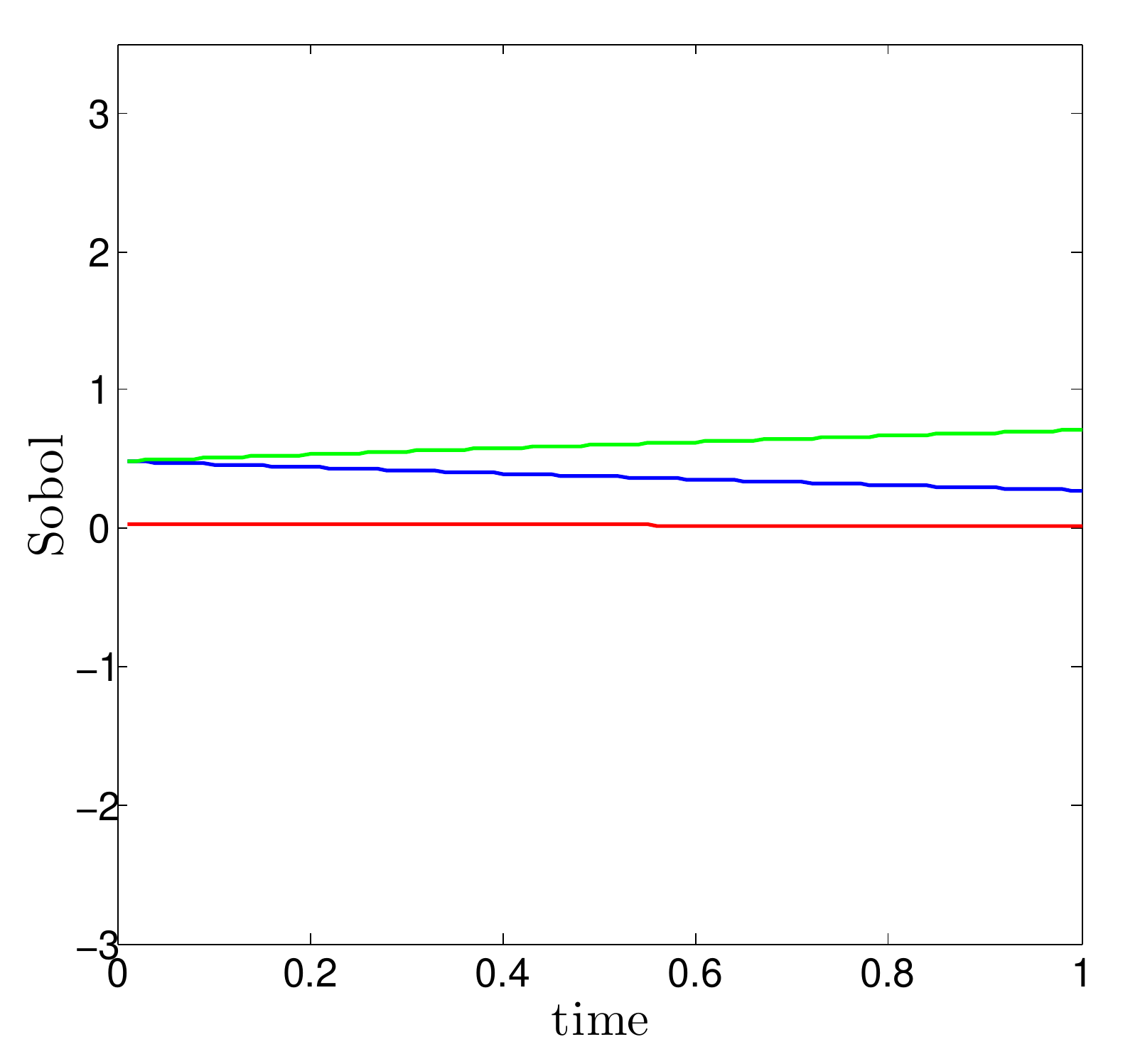}
\includegraphics[scale=0.21]{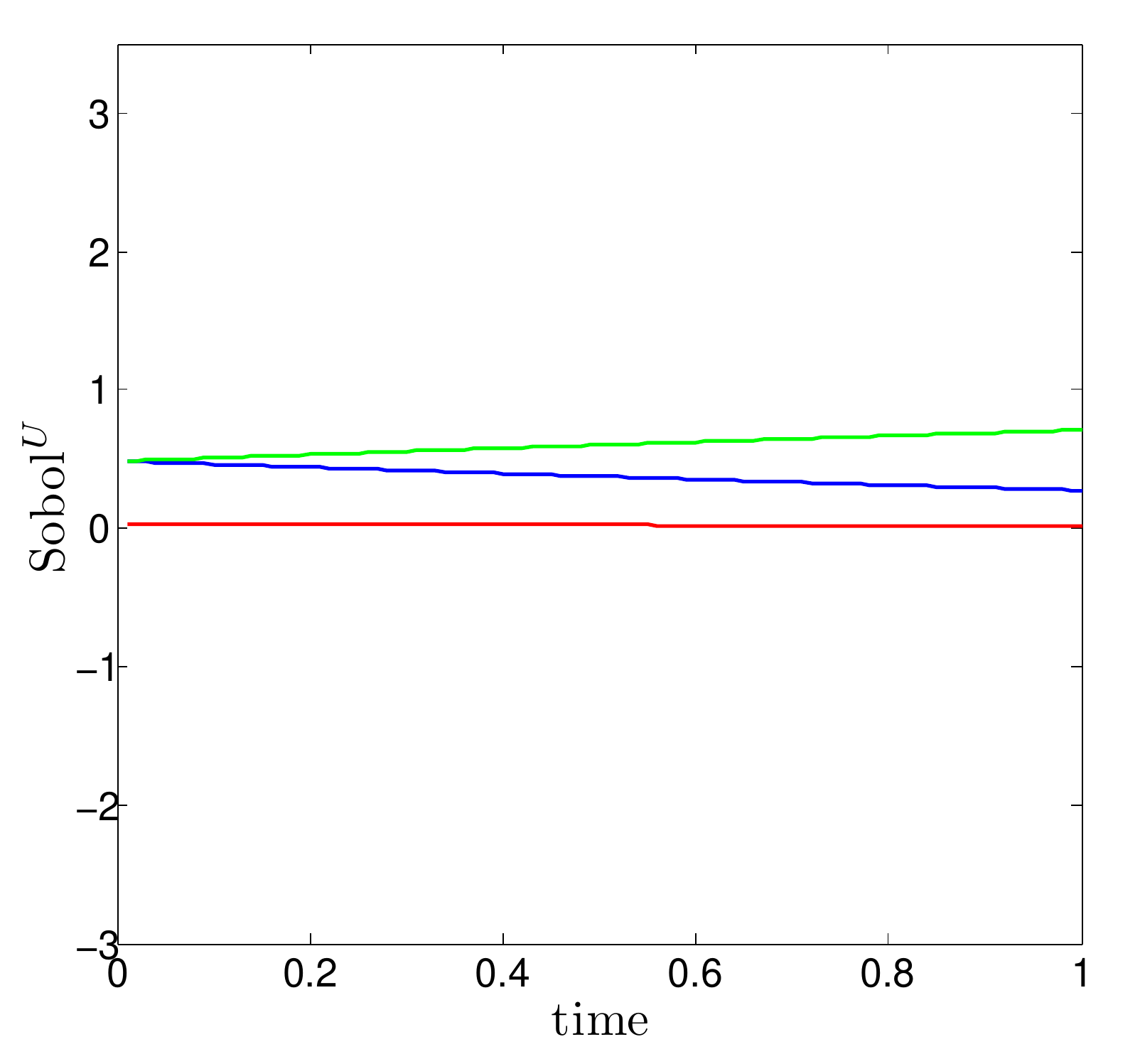}
\includegraphics[scale=0.21]{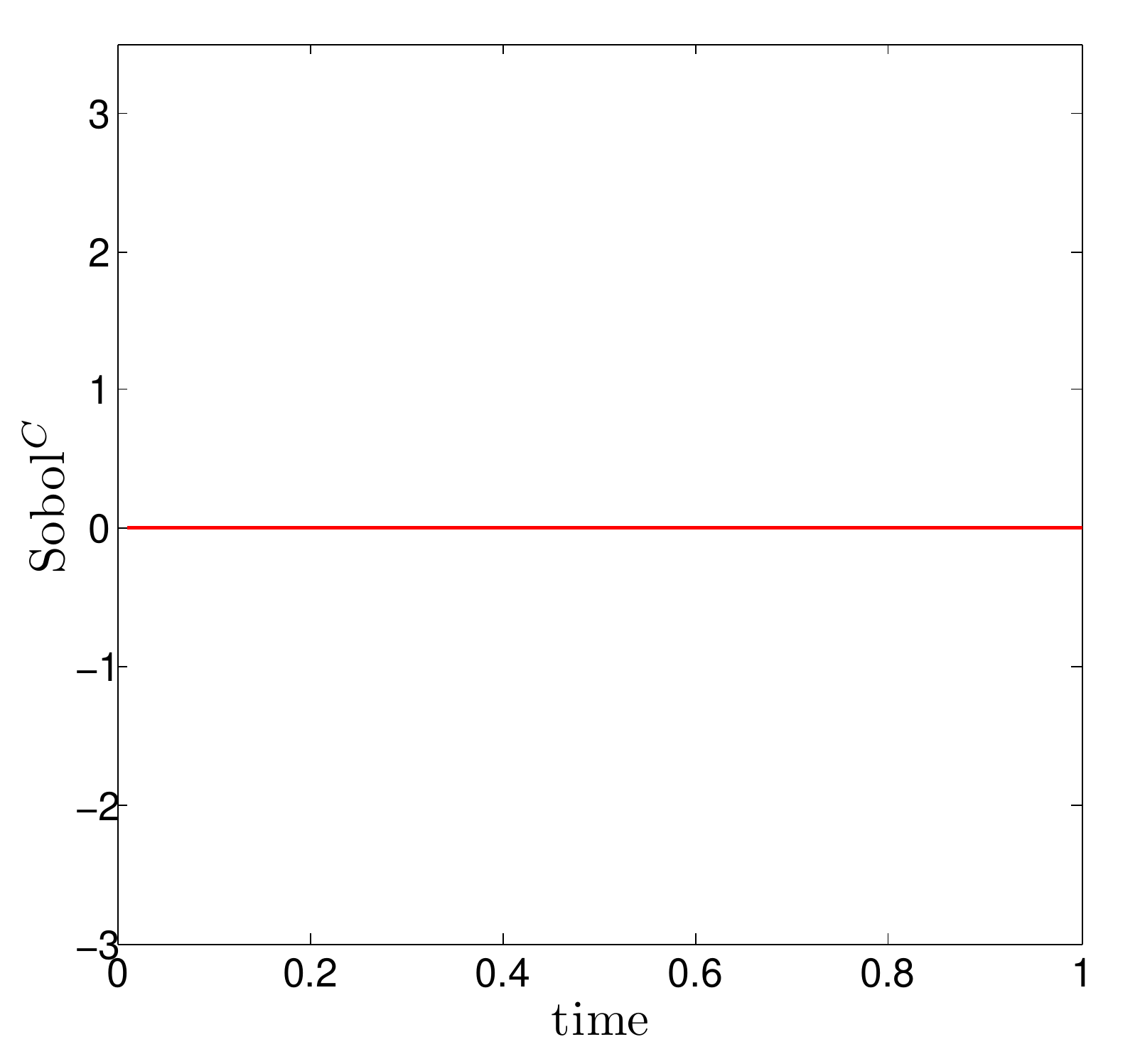}\parbox{2cm}{\vspace{-3cm} $\varrho=0$}\\
\includegraphics[scale=0.21]{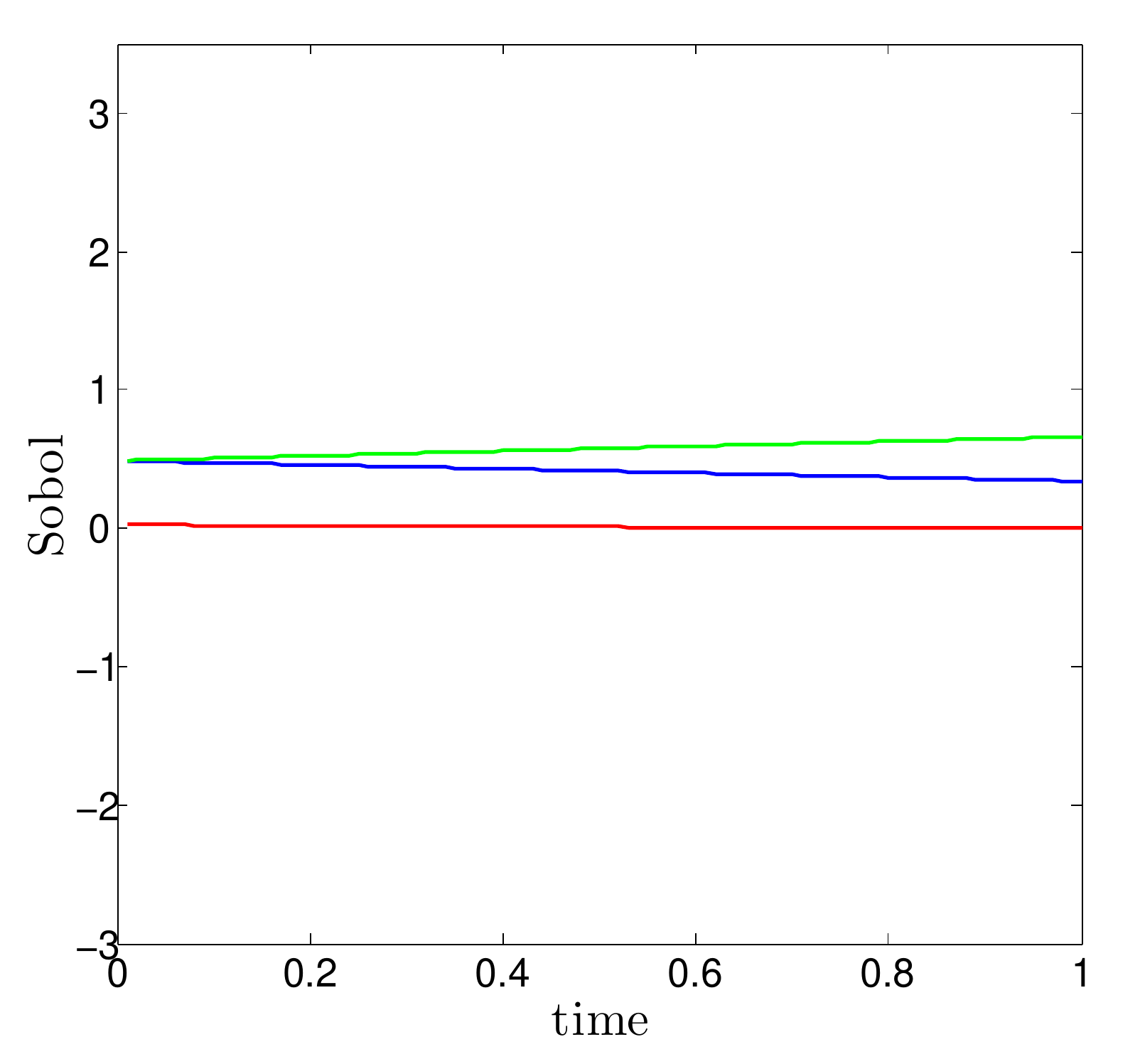}
\includegraphics[scale=0.21]{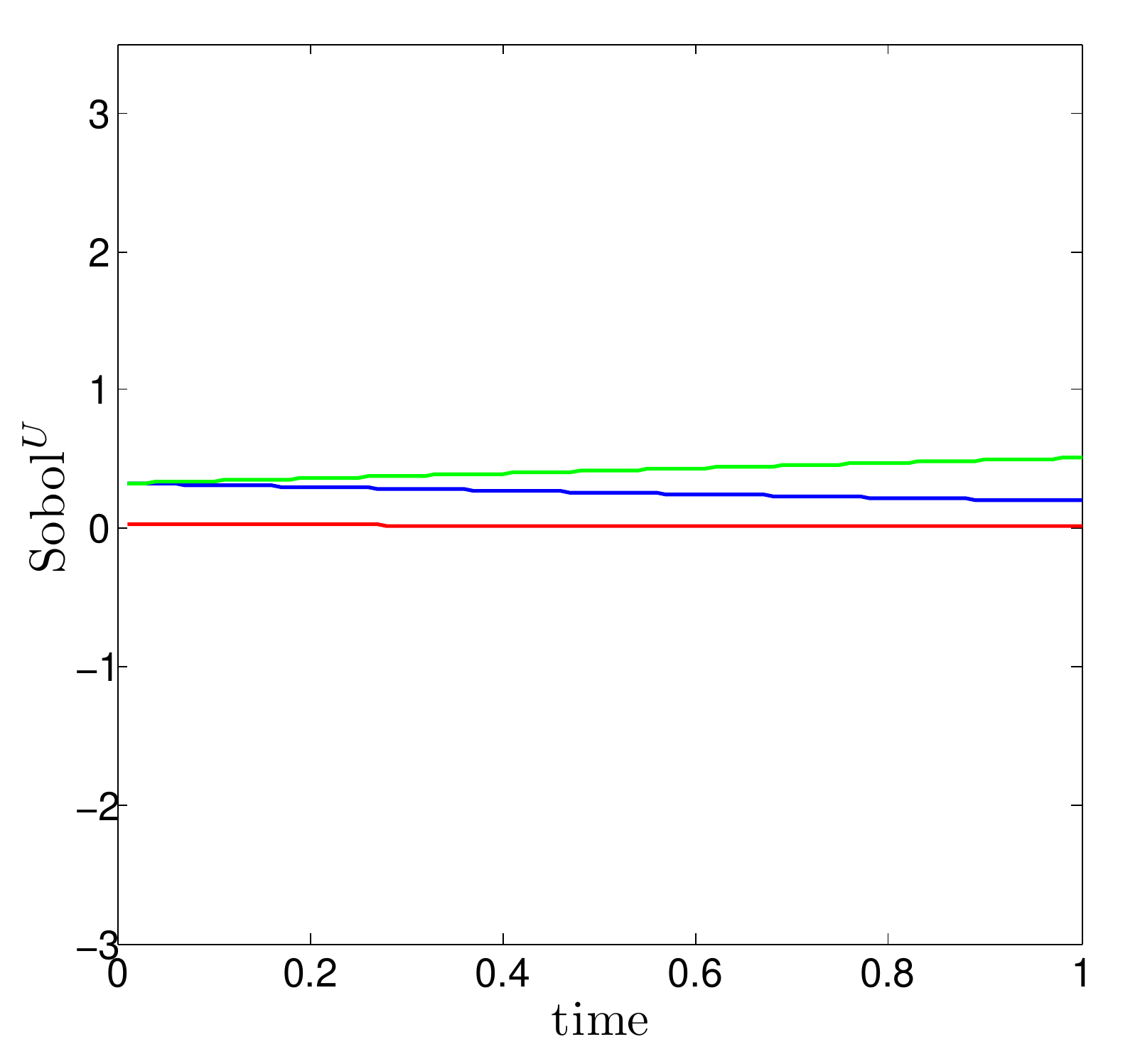}
\includegraphics[scale=0.21]{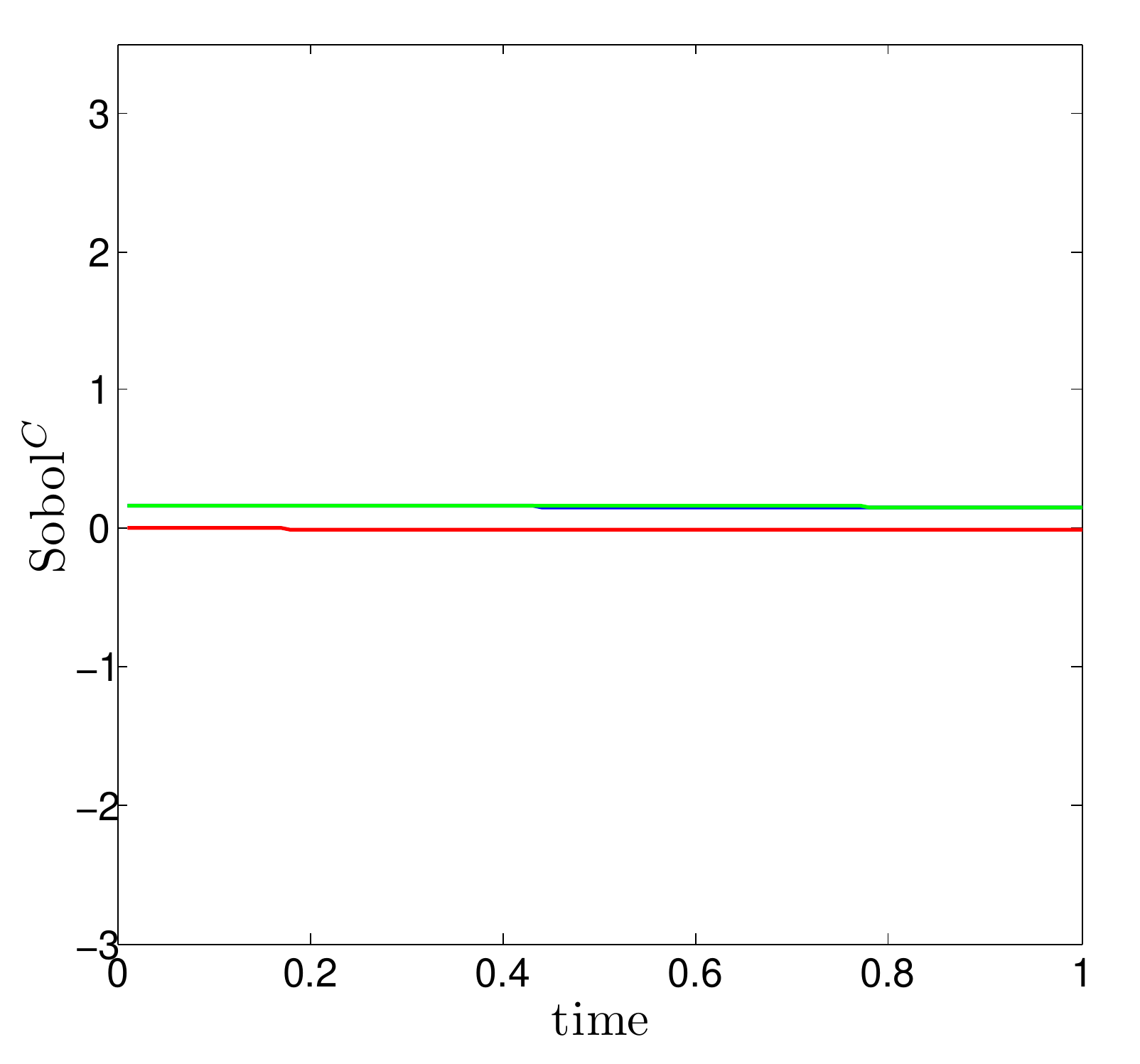}\parbox{2cm}{\vspace{-3cm} $\varrho=0.5$}\\
\includegraphics[scale=0.21]{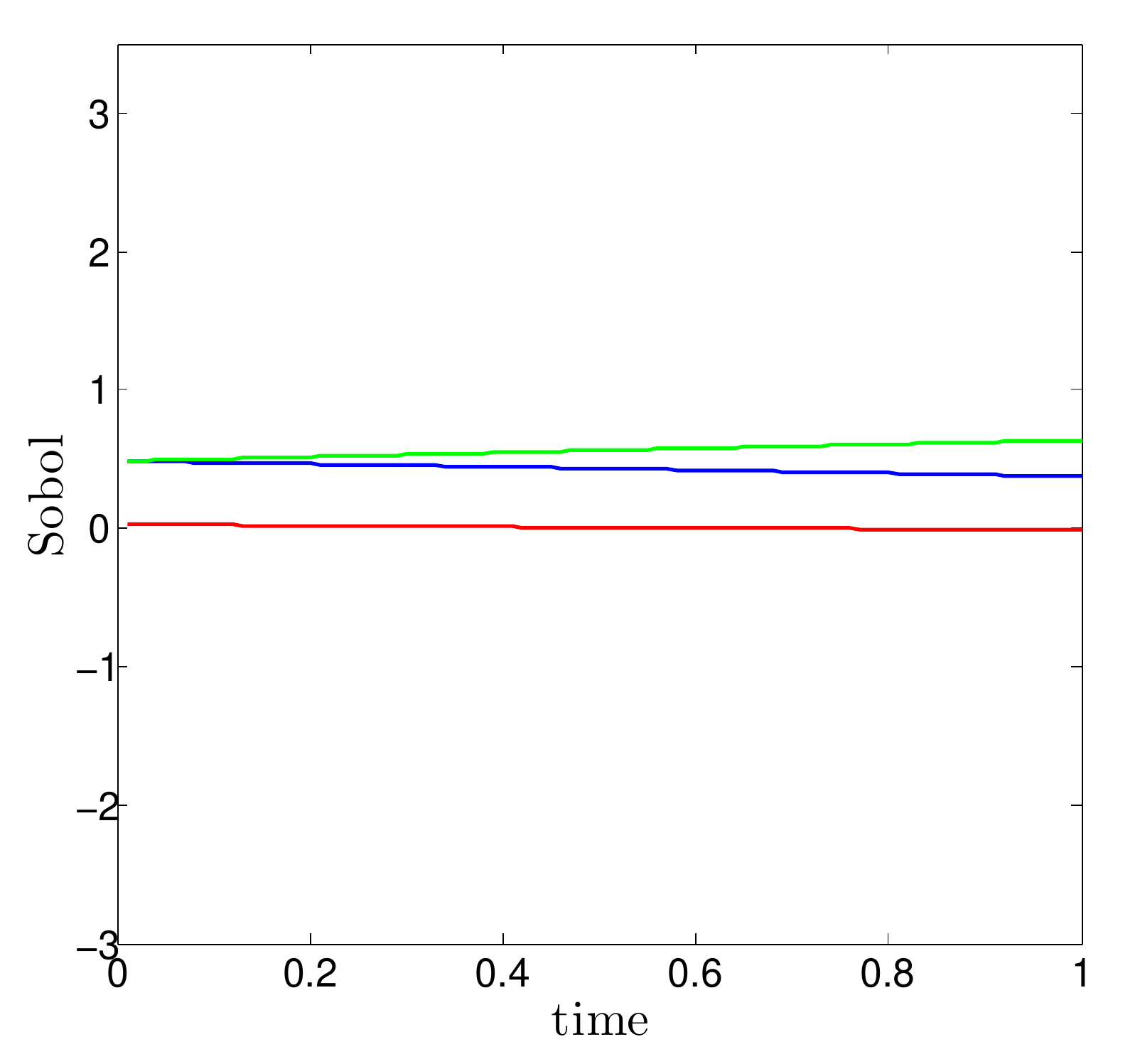}
\includegraphics[scale=0.21]{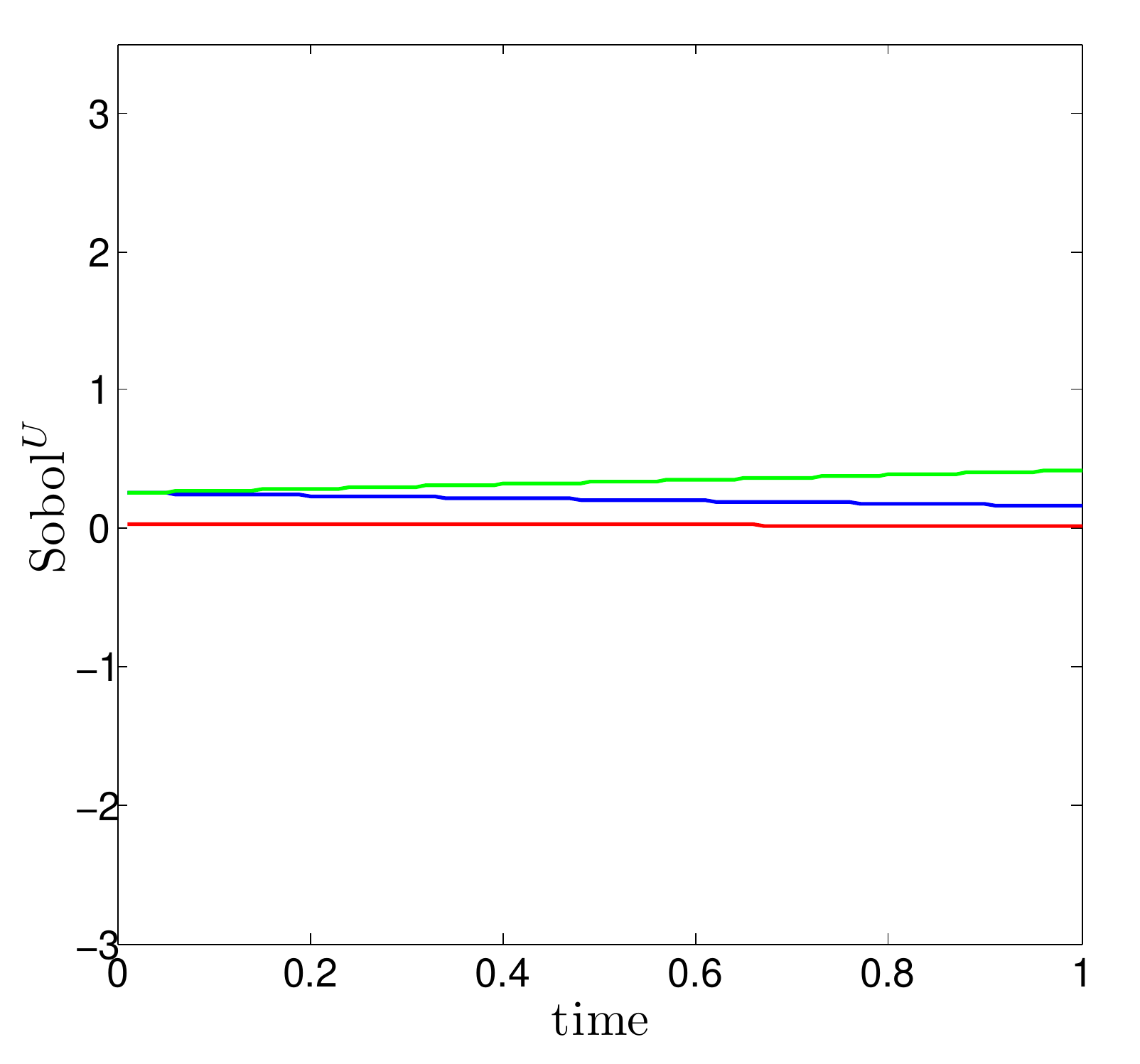}
\includegraphics[scale=0.21]{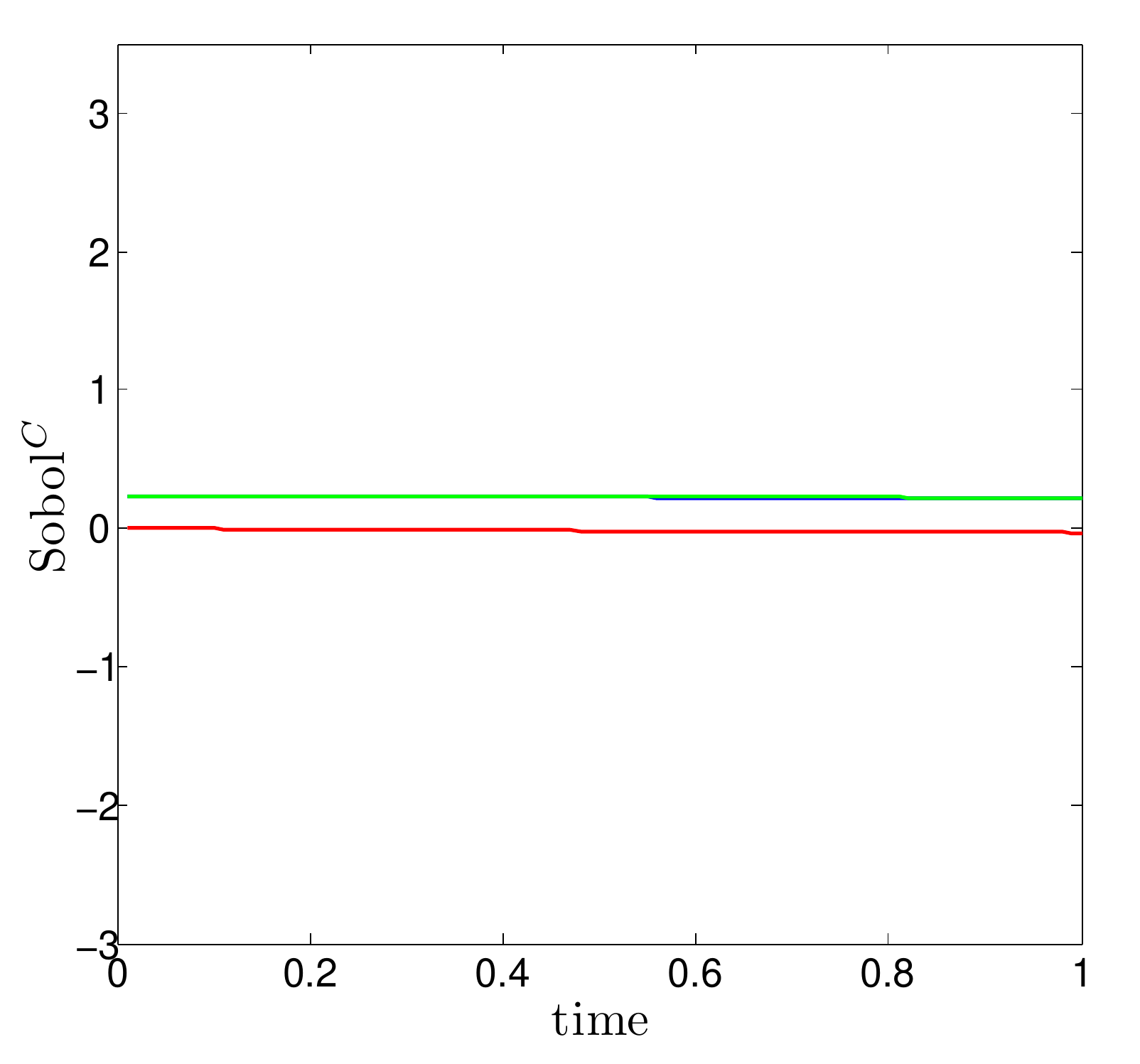}\parbox{2cm}{\vspace{-3cm} $\varrho=0.9$}\\
\end{tabular}
\caption{Dynamic plots of the Sobol' indices of decay equation (\ref{equ:Example1}). 
Left: Sobol' index, middle: contribution from the respective variables, right: contribution due to correlation with
other variables.}
\label{fig:Sobol_Ex1}
\end{figure}

Figures \ref{fig:Pol}-\ref{fig:Pol2} show all polynomials of order 0, 1 and 2 for correlation coefficients $\varrho = \{0,\pm 0.5\}$. 
The plots show clearly the non-uniqueness of the orthogonal basis, viz., the influence of the ordering  of the monic polynomials (\ref{equ:pol}): the first polynomials of a new order, $\Phi_{\{0,1,3\}}$,  are independent of the correlation coefficient and only determined by the first monic polynomial $e_{\{0,1,3\}}$ of that order, the next ones of the same order, $\Phi_{\{2,4,5\}}$, are different and dependent on the correlation.
In Figure \ref{fig:Ex1-MSd} we show the stochastic mean and standard deviation of the solution of 
Eq. (\ref{equ:Example1}) for the seven values of the correlation coefficient 
 using a PCE expansion  of 45 terms, which corresponds to a polynomial order of 8 
and the polynomial basis calculated with the proposed method. The results are plot accurate,  in the cases of no or full correlation with the ones obtained with the Hermite polynomial basis and same order, and 
with Monte Carlo simulations for the correlated distributions.
%\ERASE{ For the latter $10^5$ samples were needed to reach this accuracy [[\Maria{ I don't know if that is correct, because with $10^5$ MC points the same accuracy is not reached, we only can see how MC converges to the same solution}]].} 
The importance of propagating the {\em true} input distribution is especially seen
in the standard deviation, where neither uncorrelated nor fully correlated give a reasonable approximation
for a distribution with correlation $\pm 0.5$. This is also reflected in the Sobol' indices.
%They have been computed again for the seven values of the correlation coefficient 
%using the same PCE expansion order that was used for the mean and the standard deviation.
Figure \ref{fig:Sobol_Ex1} shows for all correlation coefficients 
the evolution of the Sobol' indices (\ref{eq:Si}) over time and the decomposition into the uncorrelated 
(\ref{eq:Si_u}) and the correlated (\ref{eq:Si_c}) part. These plots and Table~\ref{tab:Sobol_tf_Ex1} show how nicely the
contributions of the variables itself and of the correlated variables follow the strength and sign of the 
correlation coefficient $\varrho$.\\
{\bf Remark.} The fully correlated and uncorrelated results do not necessarily bound the correlated ones.
E.g., for a delta function response which is located just off the diagonal, full correlation results in a
zero variance, and no correlation gives a lower variance than correlation.

\begin{table}[htbp!]
\centering
\begin{tabular}{|r|l|ccc|}
\hline\hline
 $\varrho$& & $S_i$ & $S_i^u$ & $S_i^c$ \\
\hline
    &$S_1$&              -0.337996430555719&   1.079010380744548 &   -1.417006811300311\\
-0.9&$S_2$&                1.169072254592321 &   2.794020096839133 &    -1.624947842246789\\
    &$S_{12}$&     0.168924175966794 &    0.134133434465047&    0.034790741501840\\
\hline
&Sum&                             1.000000000003396&   4.007163912048728&  -3.007163912045260\\   
\hline
    &$S_1$&               0.125070219450471&    0.462450913298287&  -0.337380693847799\\
-0.5&$S_2$&                0.819760985060438&   1.197491238203163&   -0.377730253142723\\
    &$S_{12}$&    0.055168795489075&    0.037711983172792&    0.017456812316329\\
\hline
&Sum&     0.999999999999984&  1.697654134674242&   -0.697654134674193\\
\hline
    &$S_1$&             0.273826682991019 & 0.273826682991019&  0.000000000000004 \\
  0 &$S_2$&              0.709059149322051 &  0.709059149322055& -0.000000000000000\\
    &$S_{12}$&   0.017114167686936&   0.017114167686939&  0.000000000000015\\
\hline
&Sum&     1.000000000000006&   1.000000000000013&   0.000000000000019\\   
\hline
    &$S_1$&             0.340423008768462&   0.196827595099748 &  0.143595413668727\\
 0.5&$S_2$&              0.660699595355148 &  0.509674242072905&   0.151025353282240 \\
    &$S_{12}$&  -0.001122604123609 &  0.016050907606319&  -0.017173511729914\\
\hline
&Sum&    1.000000000000001&  0.722552744778972&  0.311794278680881\\   
\hline
    &$S_1$&              0.374323408284474 &   0.161822550653440 &   0.212500857631028\\
 0.9&$S_2$&                0.636757390356394 &   0.419027897208577&   0.217729493147815\\
    &$S_{12}$&    -0.011080798640891&  0.020103679953051&    -0.031184478593928\\
\hline
&Sum&    0.999999999999977&  0.600954127815068&   0.399045872184915\\   
\hline\hline
\end{tabular}
\caption{Sobol' indices and the contributions from variables itself and from correlated variables in the
final timepoint $t_f=1$. The respective contributions follow nicely the strength and sign of the correlation
coefficient.}
\label{tab:Sobol_tf_Ex1}
\end{table}

To answer the question whether the method has the same favourable convergence behavior as the original expansion
for independent normal random variables (cf. \cite{CandM}),
we have performed an error convergence study using the following error measures
\begin{equation}
\varepsilon_{\mu}(p)=\left|\mu_{p}(t_{f})-\mu_{exact}(t_{f})\right|, \qquad   
\varepsilon_{\sigma}(p)=\left|\sigma_{p}(t_{f})-\sigma_{exact}(t_{f})\right|,
\end{equation}
where $t_{f}=1$ is the final time point,
$\mu_{p}$  and  $\sigma_{p}$  correspond to the mean and the standard deviation obtained with 
PC expansion of order $p$, and 
$\mu_{exact}$  and $\sigma_{exact}$ are the reference solutions, in this case obtained with $p=8$. 

\begin{figure}[ht!]
   \centering
   \subfloat[Mean]{
        \includegraphics[width=6.5cm,height=6.5cm]{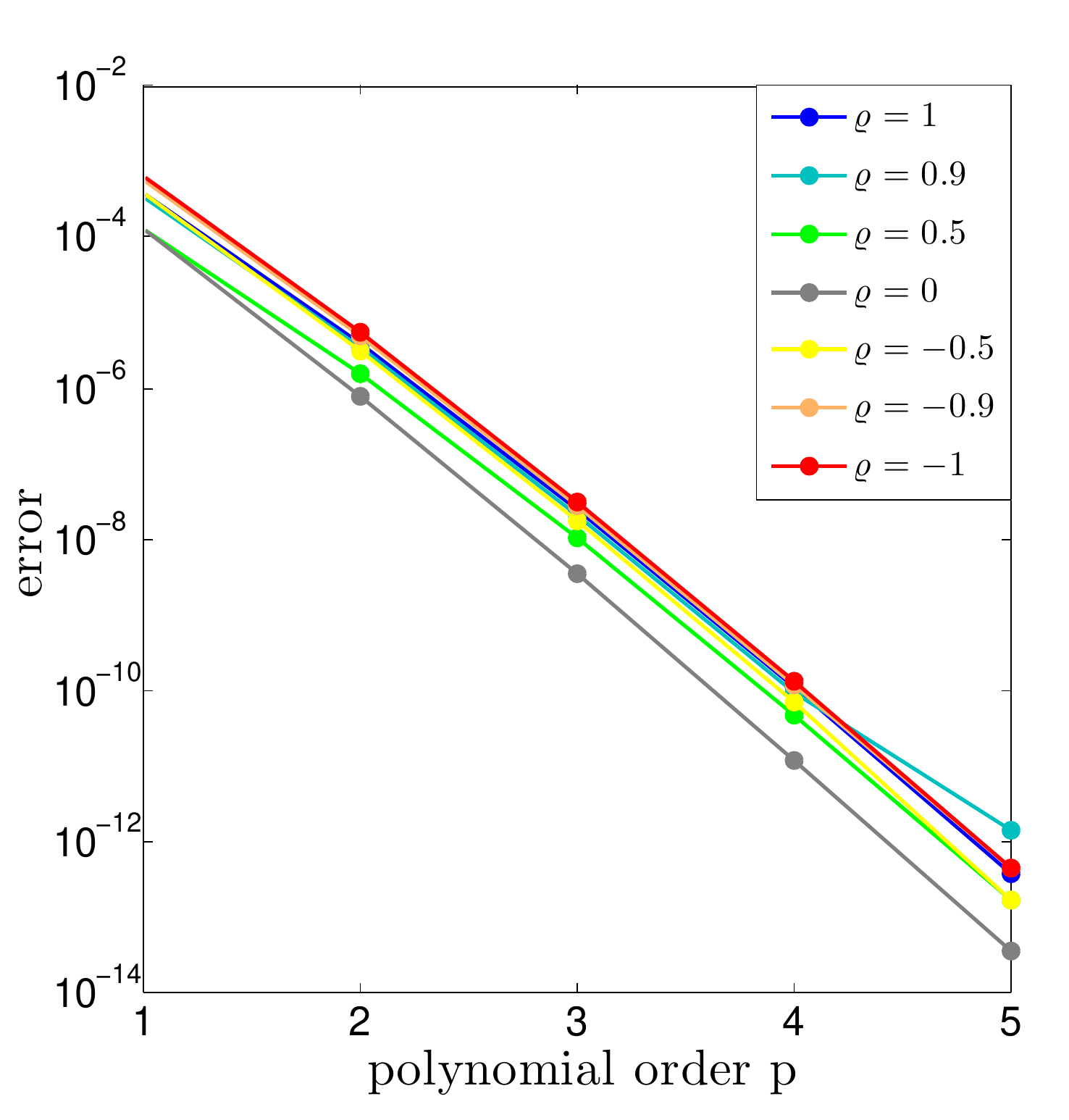}}
   \hspace{0.01\linewidth}
   \subfloat[Standard Deviation]{
        \includegraphics[width=6.5cm,height=6.5cm]{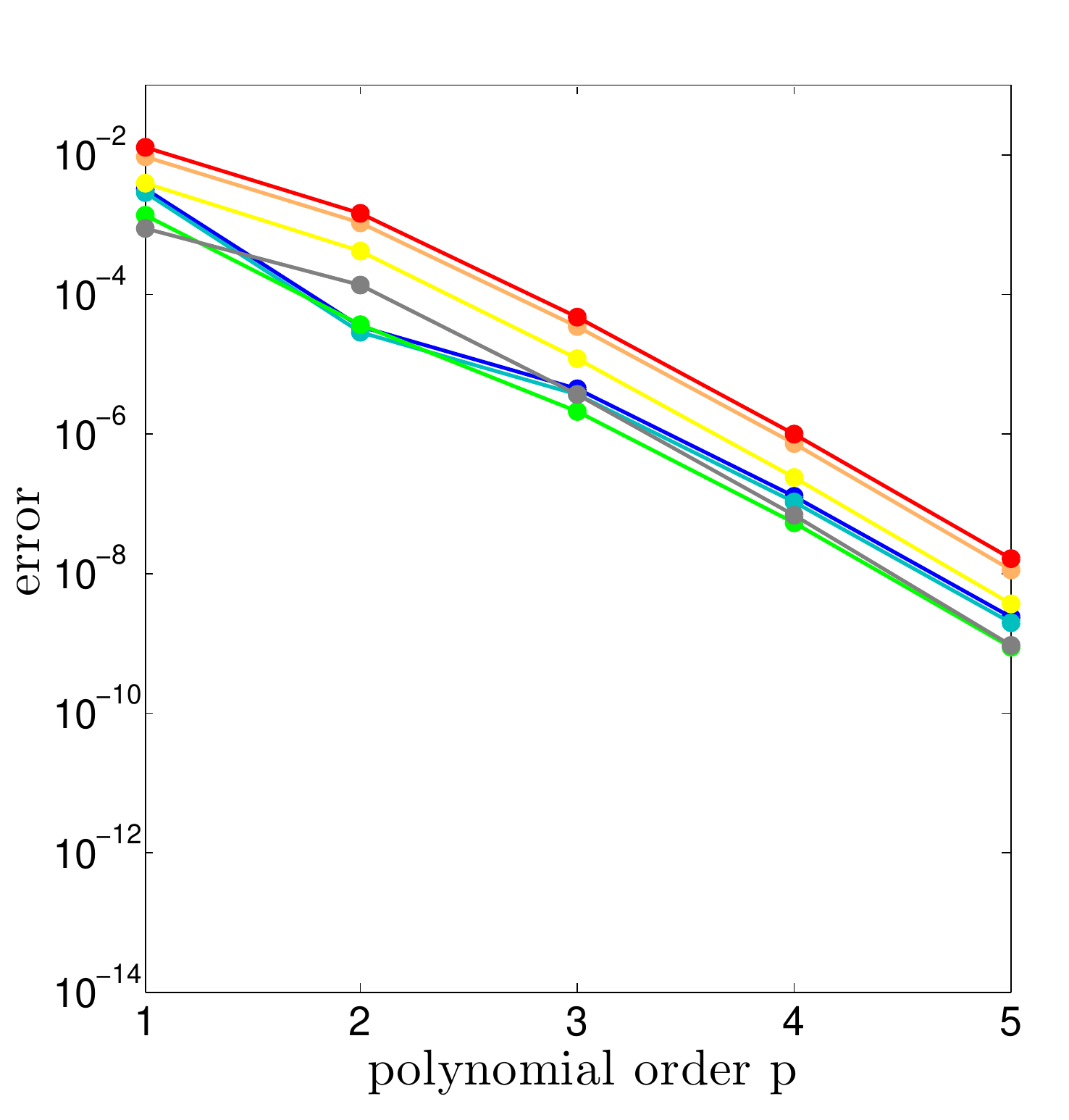}}
   \hspace{0.01\linewidth}
	\caption{Convergence of the error in mean and standard deviation for the 
	decay equation (\ref{equ:Example1})}
   \label{fig:Decay-Error-MS} 
\end{figure}
As we can see in Figure \ref{fig:Decay-Error-MS} the convergence rate of the error for increasing expansion
order is for all correlation coefficients $\varrho$ the same, showing that  the convergence rate of 
the original PCE method for independently normally distributed input variables is conserved. 

\subsection{Enzymatic reaction}
Consider the canonical enzymatic reaction
\begin{eqnarray*}
\E + \S &\equil{1}{2}& \C\\
\C &\product{3}& \E +\P 
\end{eqnarray*}
where the concentrations of the substrate \S,
the enzyme \E, the complex \C, and the product \P\ are the state variables and
$k_1$, $k_2$, and $k_3$ are the kinetic rate constants.
This example has been used in a previous paper \cite{FEBS} to show the effect of noisy data on the 
reliability of the estimated parameters using linear regression. However, in most cases biologists
are not so much interested in accurate model parameters, but in an accurate solution of the QoI, in this
case the concentration of the product \P. We now can propagate the
uncertainty in the parameters, given by the Fisher information matrix, so including correlation
between the kinetic parameters, through the model to obtain the pdf for the QoI \c\P.
The mathematical model for this enzymatic reaction is given by the following system of ODEs:
\begin{eqnarray}
\frac{d\c\S(t)}{dt} &=& -k_1\c\E(t)\c\S(t)+k_2\c\C(t) \nonumber\\
\frac{d\c\C(t)}{dt} &=& k_1\c\E(t)\c\S(t)-k_2\c\C(t)-k_3\c\C(t) \nonumber\\
\frac{d\c\E(t)}{dt} &=& -k_1\c\E(t)\c\S(t)+k_2\c\C(t)+k_3\c\C(t) \\
\frac{d\c\P(t)}{dt} &=& k_3\c\C(t)\nonumber
\label{equ:FEBS}
\end{eqnarray}
Data for \c\C\ are available at regular intervals during the time interval [0,20] and 
the initial concentrations are known
$\c\S(0)=1$, $\c\C(0)=0$, $\c\E(0)=1$, and $\c\P(0)=0$.
The parameters $k_1$, $k_2$, and $k_3$ have been obtained by linear regression using the noisy data for \c\C.
The optimal value of the parameters in \cite{FEBS} was  $\bsk=(0.683, 0.312, 0.212)$, with (independent)
confidence intervals $\Delta=[0.033, 0.028, 0.005]$. 
To study the propagation of the uncertainty through the model, the parameters are assumed to be random and their pdf
is built based on the results in \cite{FEBS}. Thus, we assumed that the parameters follow a uniform distribution with mean
the optimal value of the parameters $\bsmu=(0.683, 0.312, 0.212)$, standard deviation $\bssi=(0.206, 0.175, 0.031)$
based on the amplitude of the confidence intervals (with multiplication factor limited to 6.25 to avoid
entering the non-physical negative parameter space),
and $C$ the correlation matrix also obtained from the 
Fisher information matrix; as in the previous example, the fully correlated case, $C_f$, 
and the uncorrelated case, $C_u$, are also studied
\begin{eqnarray*}
&& C = \begin{pmatrix}
1 & 0.9& -0.37 \\
0.9 & 1& -0.45 \\
-0.37 & -0.45& 1 \\
\end{pmatrix}\;
C_f = \begin{pmatrix}
1 & 1& -1\\
1 & 1& -1 \\
-1 & -1& 1 \\
\end{pmatrix}\;
C_u = \begin{pmatrix}
1 & 0& 0\\
0 & 1& 0 \\
0& 0& 1 \\
\end{pmatrix}\
\end{eqnarray*}

As the analytic expression for the pdf of the joint distribution of two or more correlated 
uniform variables is unknown, we cannot use the moment-generating function to  calculate all the required moments. 
Therefore, Monte Carlo integration will be used in this example.   
To generate sampling points from a correlated multidimensional uniform  distribution we used the standard
approach of computing them from the correlated normal distribution (see, e.g., \cite{Ucorr}).
%First, a normally distributed multivariate random variable $Z$ with correlation matrix $C$
%is created, $Z \sim \mathcal{N}(0,C)$  
%The multivariate
%correlated uniform random variable, $U$, is defined from the normal one by applying the transformation $U=F(Z)$, where $F$ is the cumulative
%distribution function of $Z$. Note, that the Pearson linear correlation coefficient is not necessarily invariant to transformation (only to certain linear ones), so to preserve the correlation we will
%use the Spearman correlation coefficient to create $Z$ which is invariant for any transformation. Therefore, when the correlations between variables in the problem are given by the Pearson correlation coefficient, we will have to adjust the correlations before generating the normal distribution. The exact
%relationship between Spearman correlation and Pearson correlation is  known \cite{Spearman} and given by:
%\begin{equation}
%C^{adj}=2 \sin(\frac{\pi}{6} C).
%\end{equation}
%\JGB{[Why this stuff at the end? Do we know that $C$ from Fisher is Pearson? If so, why not describe the procedure with
%$C^{adj}$ from the beginning? Might be better to just mention the ref (and I have to understand this)]}.
%Using $U$ and Monte Carlo sampling, we generate all required moments and the orthogonal polynomial base.
Again, the problem collapses to one-dimensional in the fully correlated case, correlation matrix  $C_f$.
%it is given by the following Equation: 
%\begin{eqnarray*}
%&& \frac{d\c\S(t)}{dt}=-\frac{k_{1}+0.167}{0.9765}\c\E(t)\c\S(t)+\frac{k_{1}+0.0331}{1.7088}\c\C(t) \nonumber\\
%&& \frac{d\c\C(t)}{dt}=\frac{k_{1}+0.167}{0.9765}\c\E(t)\c\S(t)-\frac{k_{1}+0.0331}
%{1.7088}\c\C(t)-\frac{(-%k_{1})+1+7.1552}{36.1010}\c\C(t) \nonumber\\
%&& \frac{d\c\E(t)}{dt}=-\frac{k_{1}+0.167}{0.9765}\c\E(t)\c\S(t)+\frac{k_{1}+0.0331}{1.7088}\c\C(t) + \frac{(-
%k_{1})+1+7.1552}{36.1010}\c\C(t)\nonumber\\
%&& \frac{d\c\P(t)}{dt}=\frac{(-k_{1})+1+7.1552}{36.1010}\c\C(t)
%\label{equ:FEBS1d}
%\end{eqnarray*}
The only parameter in this case follows the standard uniform distribution $\mathcal{U}(0,1)$.

\begin{figure}[ht!]
%   \centering
   %%----primera subfigura----
   \subfloat[\c{S}]{
        %\label{fig:Points:a}         %% Etiqueta para la primera subfigura
        \includegraphics[width=0.33\linewidth]{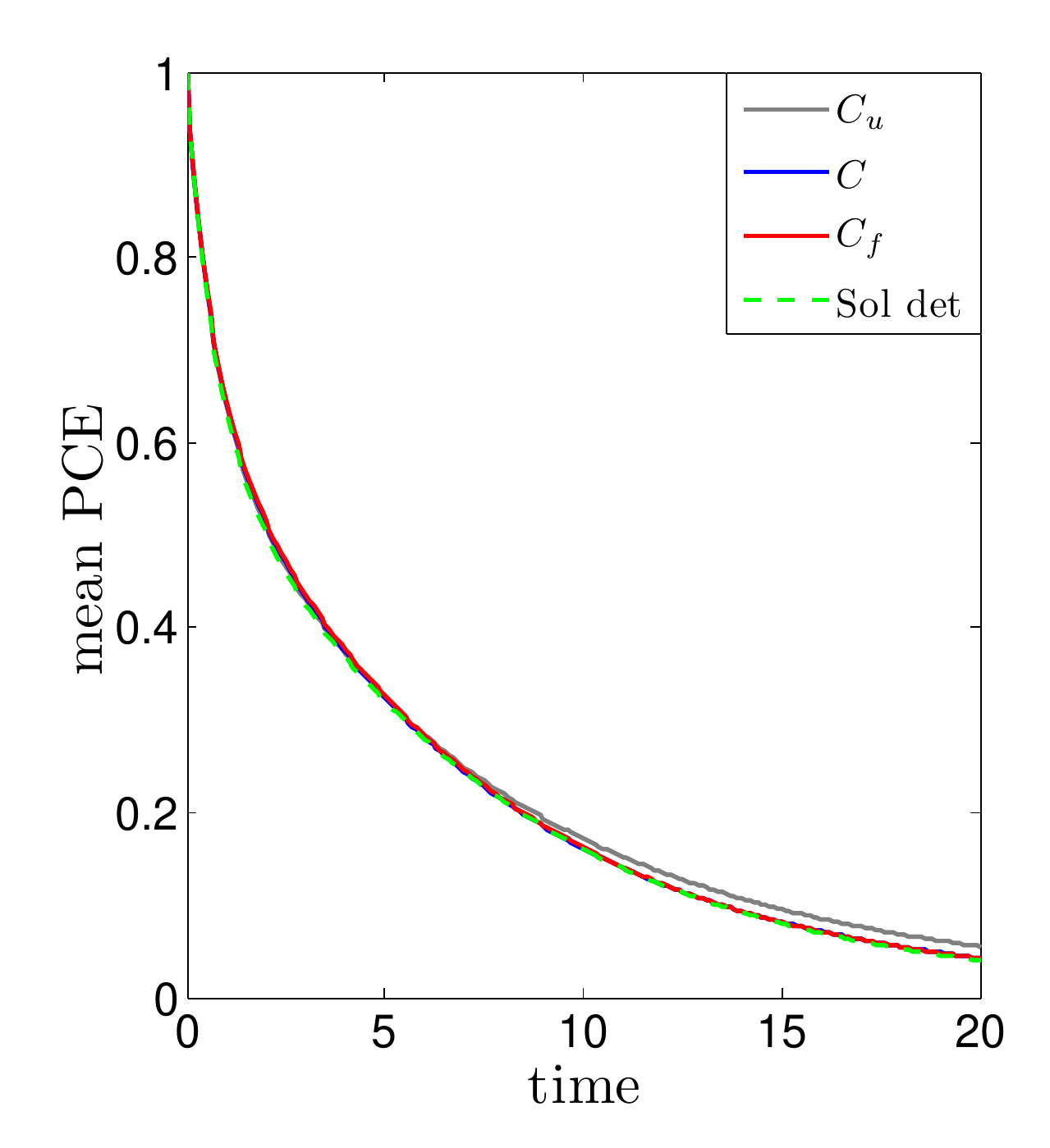}}
%   \hspace{0.01\linewidth}
   %%----segunda subfigura----
%   \subfloat[Complex C]{
        %\label{fig:Points:b}         %% Etiqueta para la segunda subfigura
%        \includegraphics[scale=0.4]{meansC}}
%   \hspace{0.01\linewidth}\\  
   \subfloat[\c{E}]{
%        %\label{fig:Points:d}         %% Etiqueta para la primera subfigura
        \includegraphics[width=0.33\linewidth]{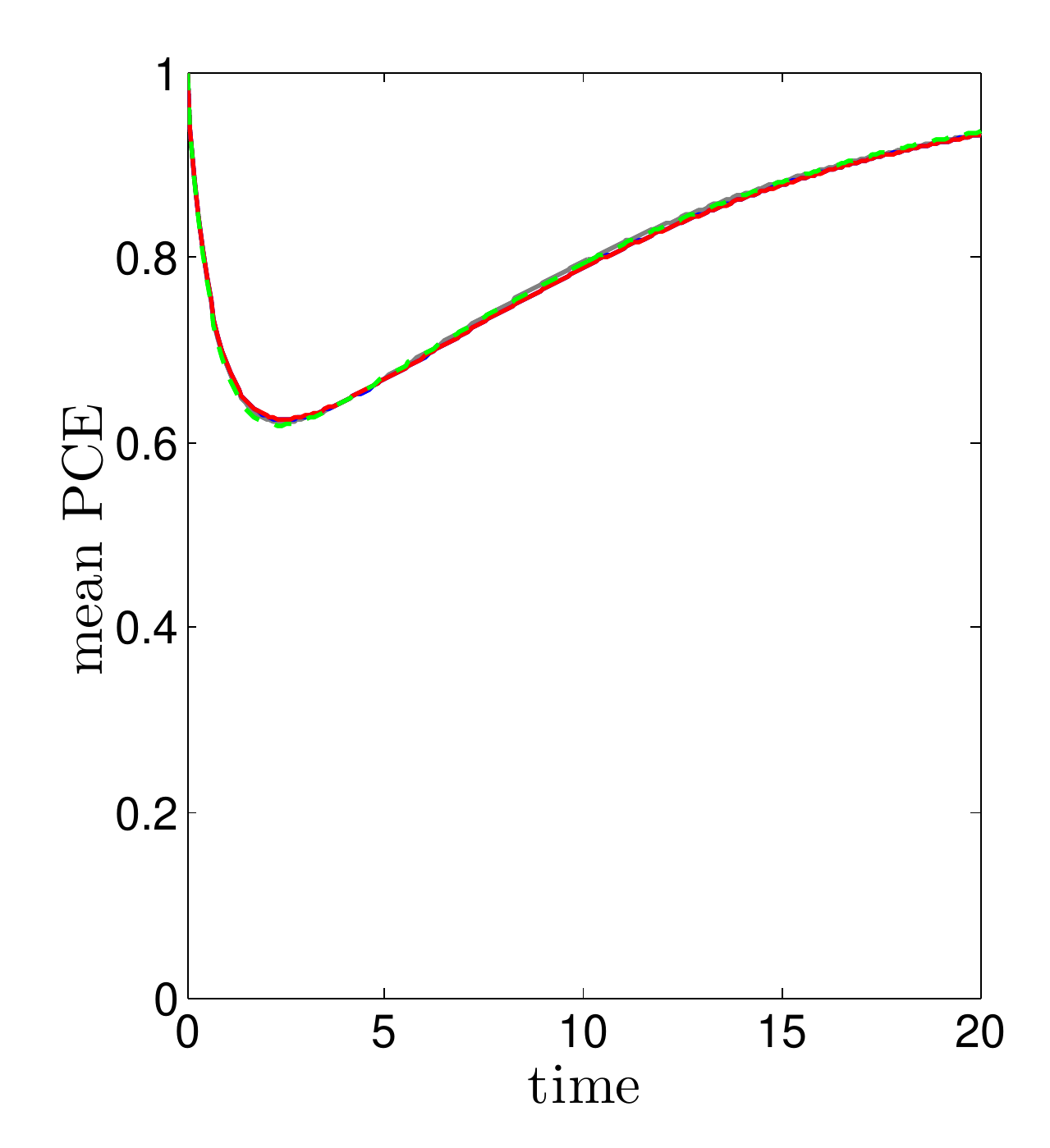}}
%   \hspace{0.01\linewidth}
   %%----segunda subfigura----
   \subfloat[\c{P}]{
       % \label{fig:Points:e}         %% Etiqueta para la segunda subfigura
        \includegraphics[width=0.33\linewidth]{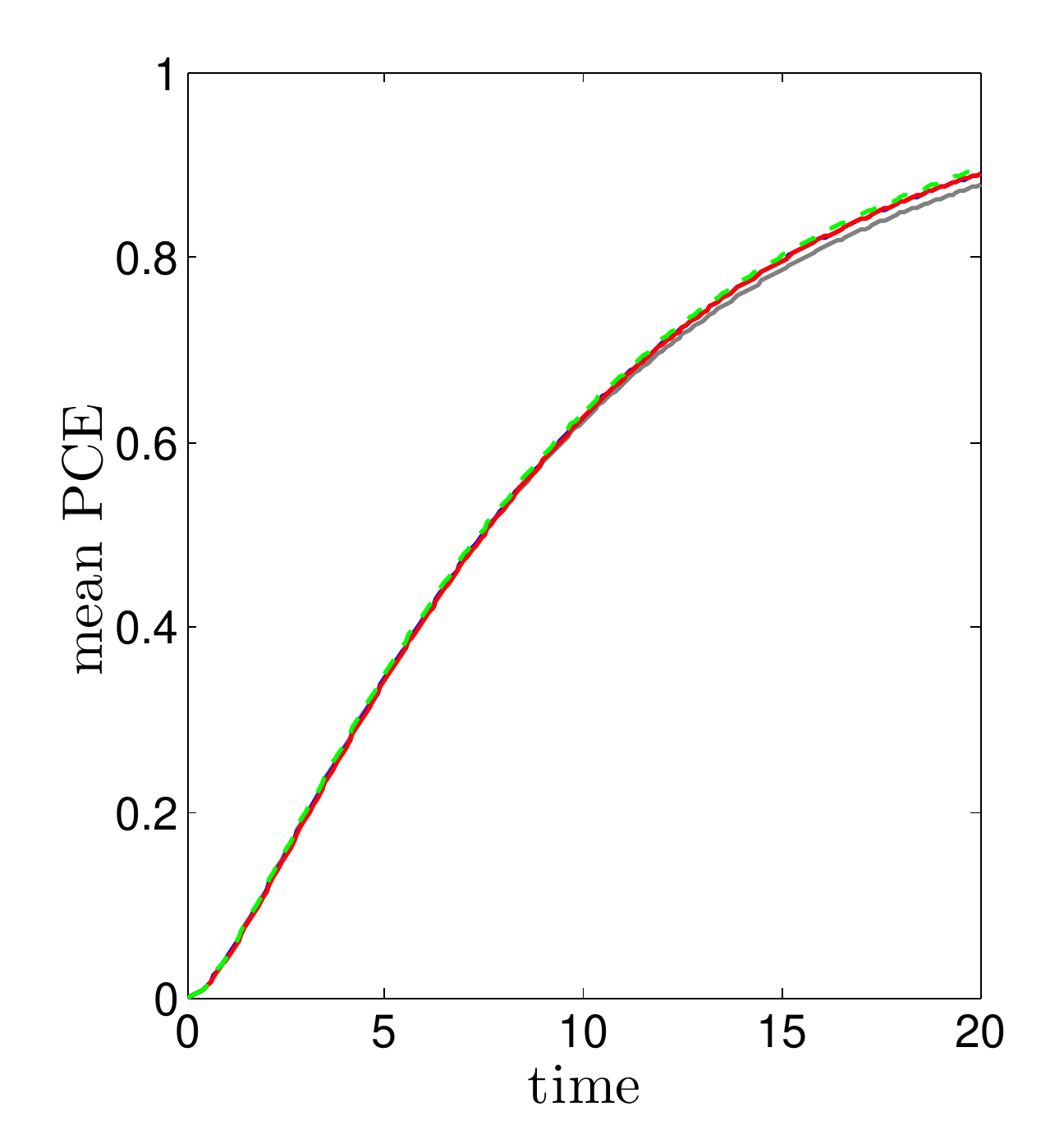}}
%				\hspace{0.01\linewidth}
   %%----segunda subfigura----
  \caption{Mean of the stochastic state variables of Eq.~(\ref{equ:FEBS}).}  
     \label{fig:Means} 
\end{figure}

\begin{figure}[ht!]
%   \centering
   %%----primera subfigura----
   \subfloat[\c{S}]{
        %\label{fig:Points:a}         %% Etiqueta para la primera subfigura
       \includegraphics[width=0.33\linewidth]{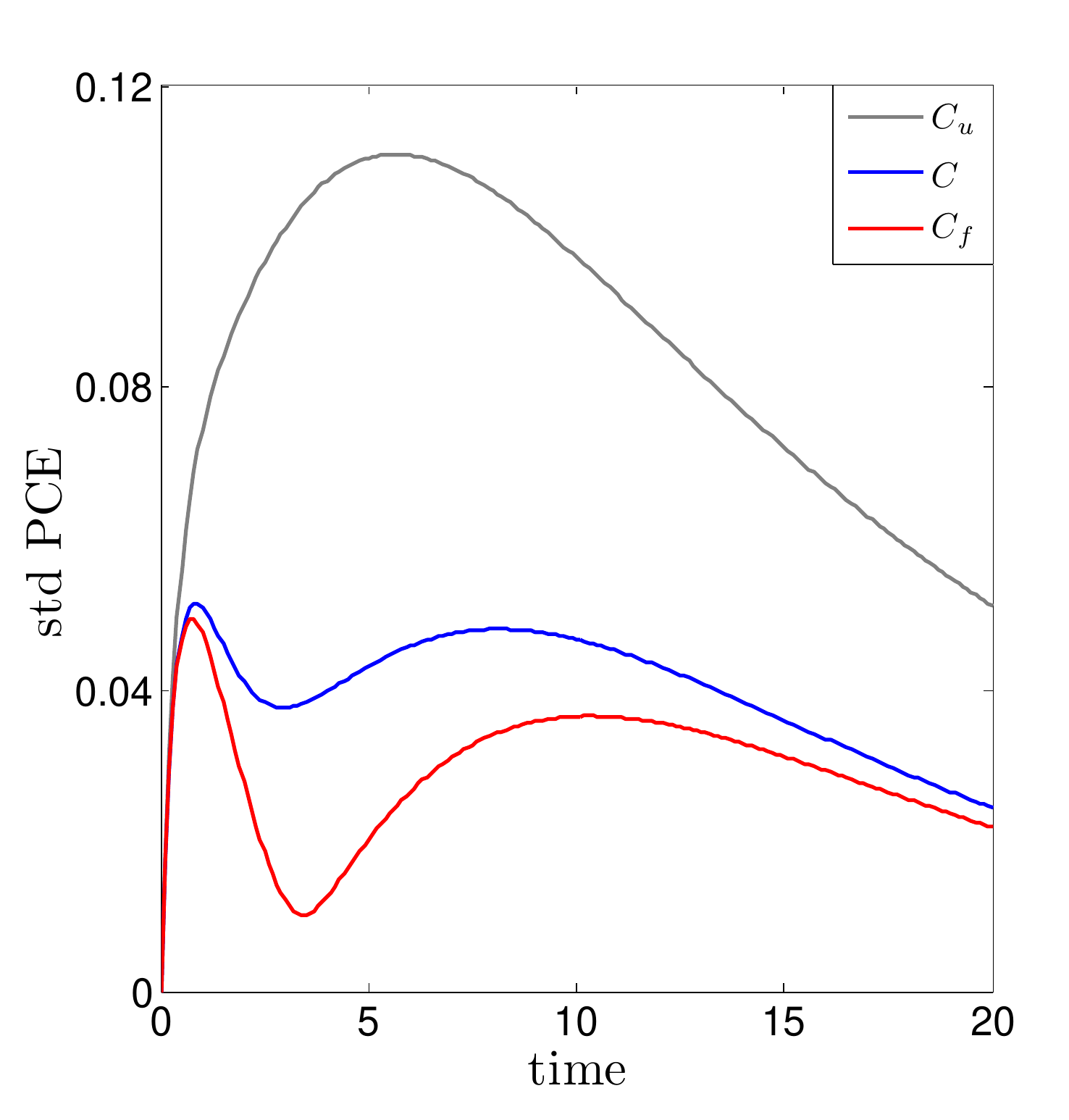}}
%   \hspace{0.01\linewidth}
   %%----segunda subfigura----
   \subfloat[\c{C} and \c{E}]{
        %\label{fig:Points:b}         %% Etiqueta para la segunda subfigura
        \includegraphics[width=0.33\linewidth]{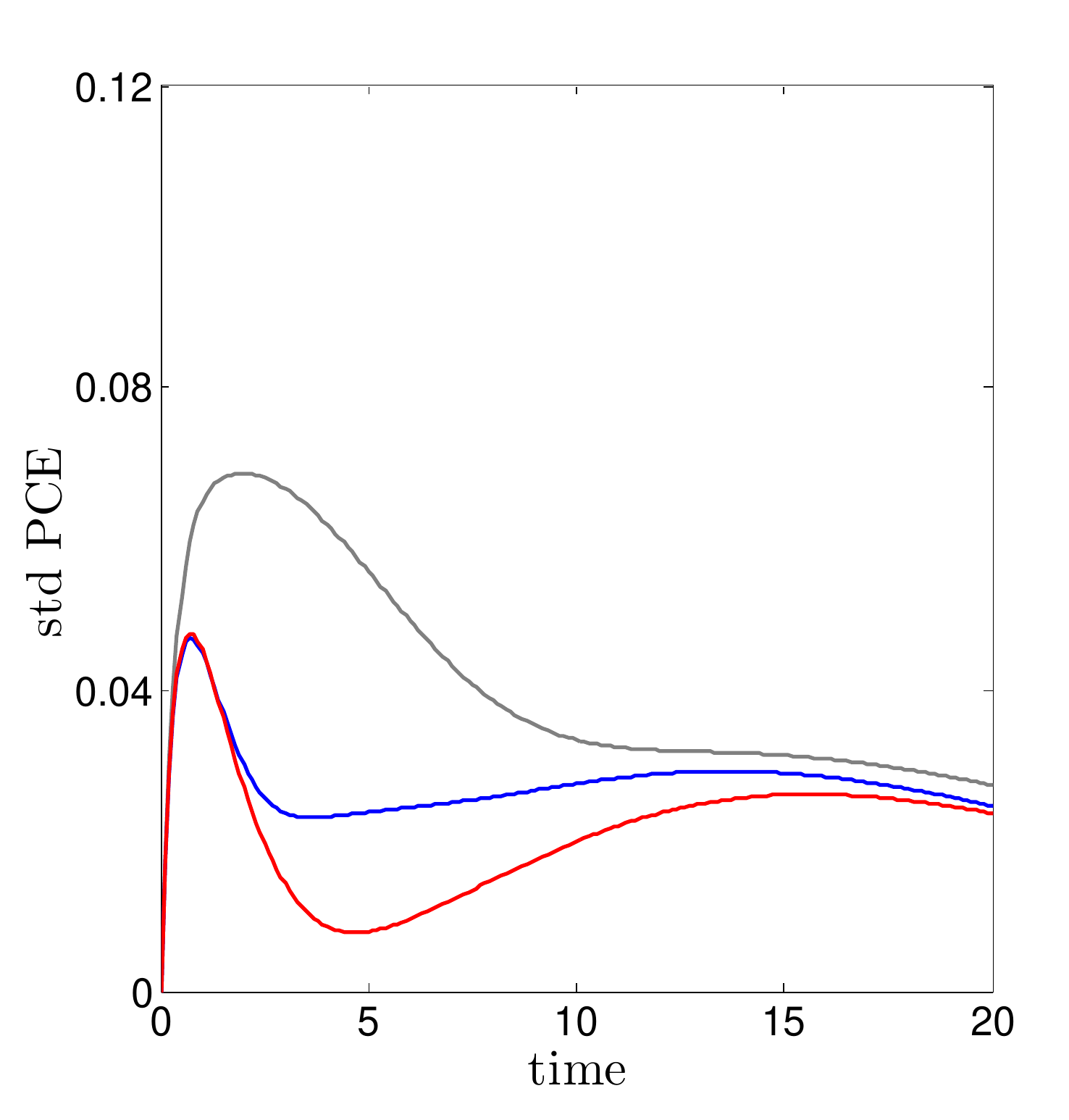}}
%   \hspace{0.01\linewidth}\\   
%   \subfloat[Enzyme E]{
       % \label{fig:Points:d}         %% Etiqueta para la primera subfigura
%      \includegraphics[scale=0.35]{stdE}}
%   \hspace{0.01\linewidth}
   %%----segunda subfigura----
   \subfloat[\c{P}]{
       % \label{fig:Points:e}         %% Etiqueta para la segunda subfigura
        \includegraphics[width=0.33\linewidth]{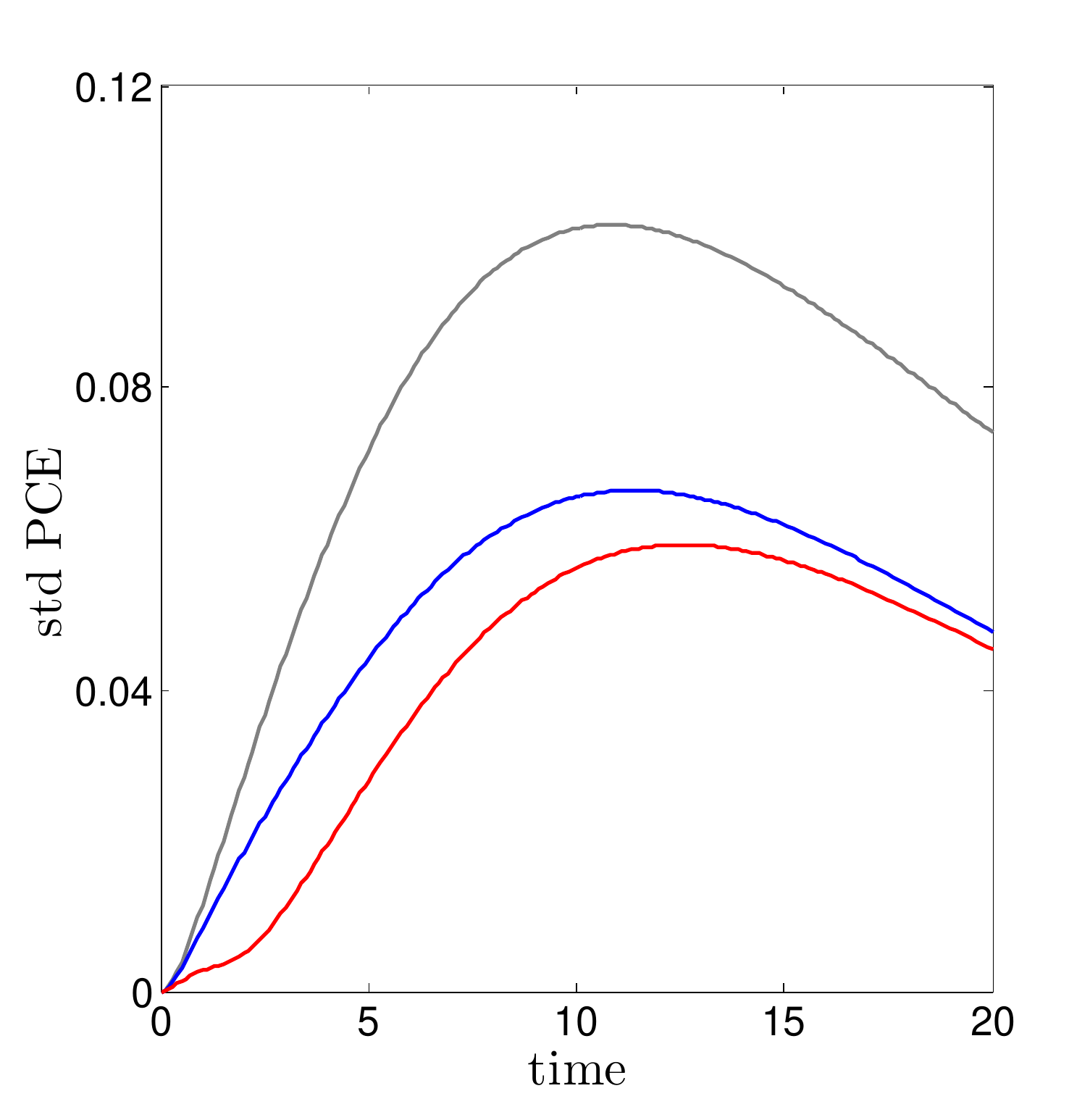}}
				\hspace{0.01\linewidth}
   %%----segunda subfigura----
 \caption{Standard deviation of the stochastic state variables of Eq.~(\ref{equ:FEBS}) 
 for the uncorrelated case (grey), the fully correlated case (red), and the correlation matrix resulting from
 the experimental errors (blue).}
   \label{fig:Stds} 
\end{figure}

For all three cases - correlated, uncorrelated, and fully correlated - we applied PCE 
with Galerkin projection to Eq. (\ref{equ:FEBS}). 
\begin{figure}[htbp!]
\begin{tabular}{ccc}
\centerline{\includegraphics[width=0.7\linewidth]{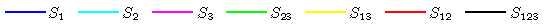}}\\
\hspace{-0.7cm}
\subfloat[\c{S}]{
\includegraphics[width=0.33\linewidth]{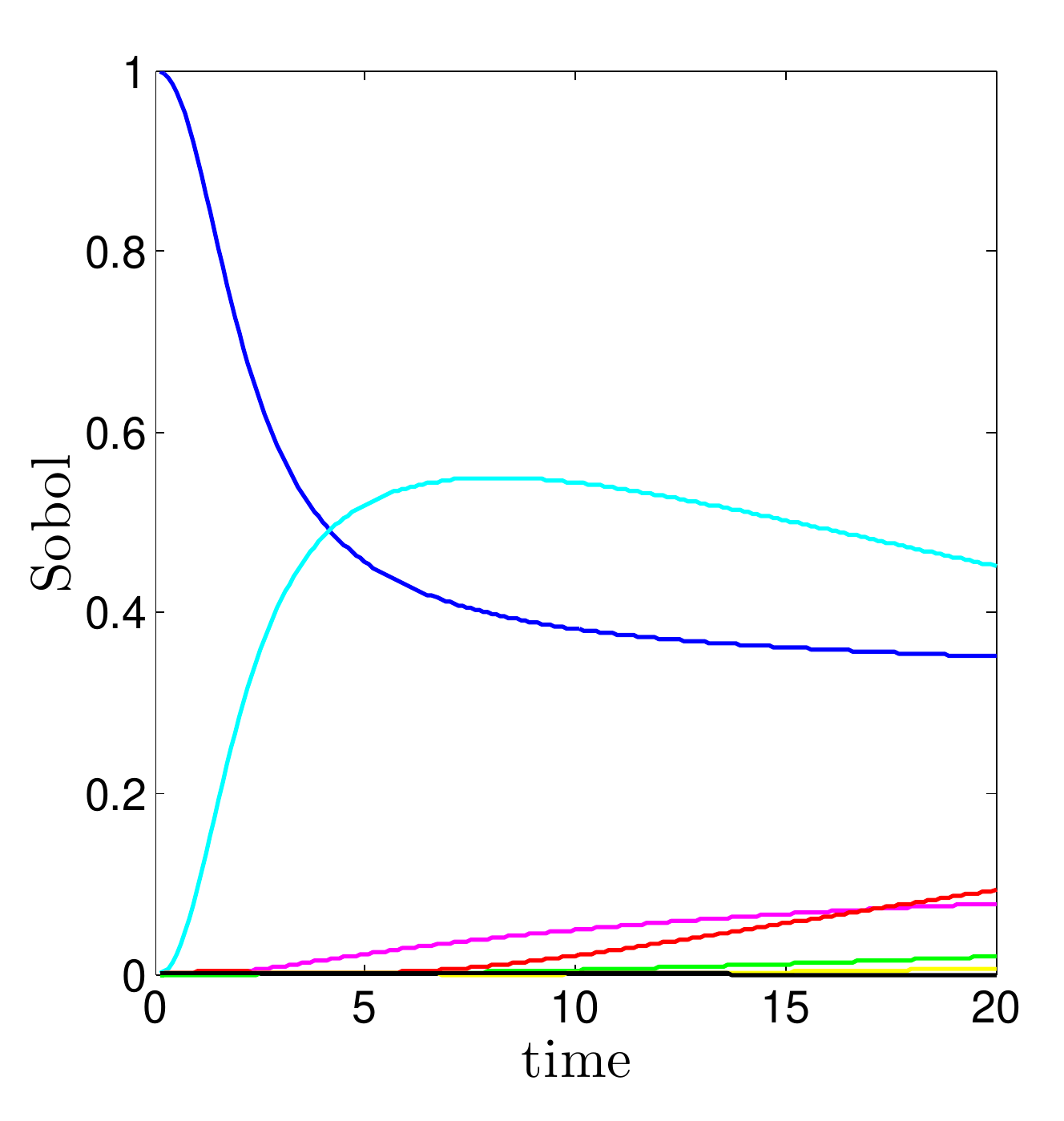}}
\subfloat[\c{C} and \c{E}]{
\includegraphics[width=0.33\linewidth]{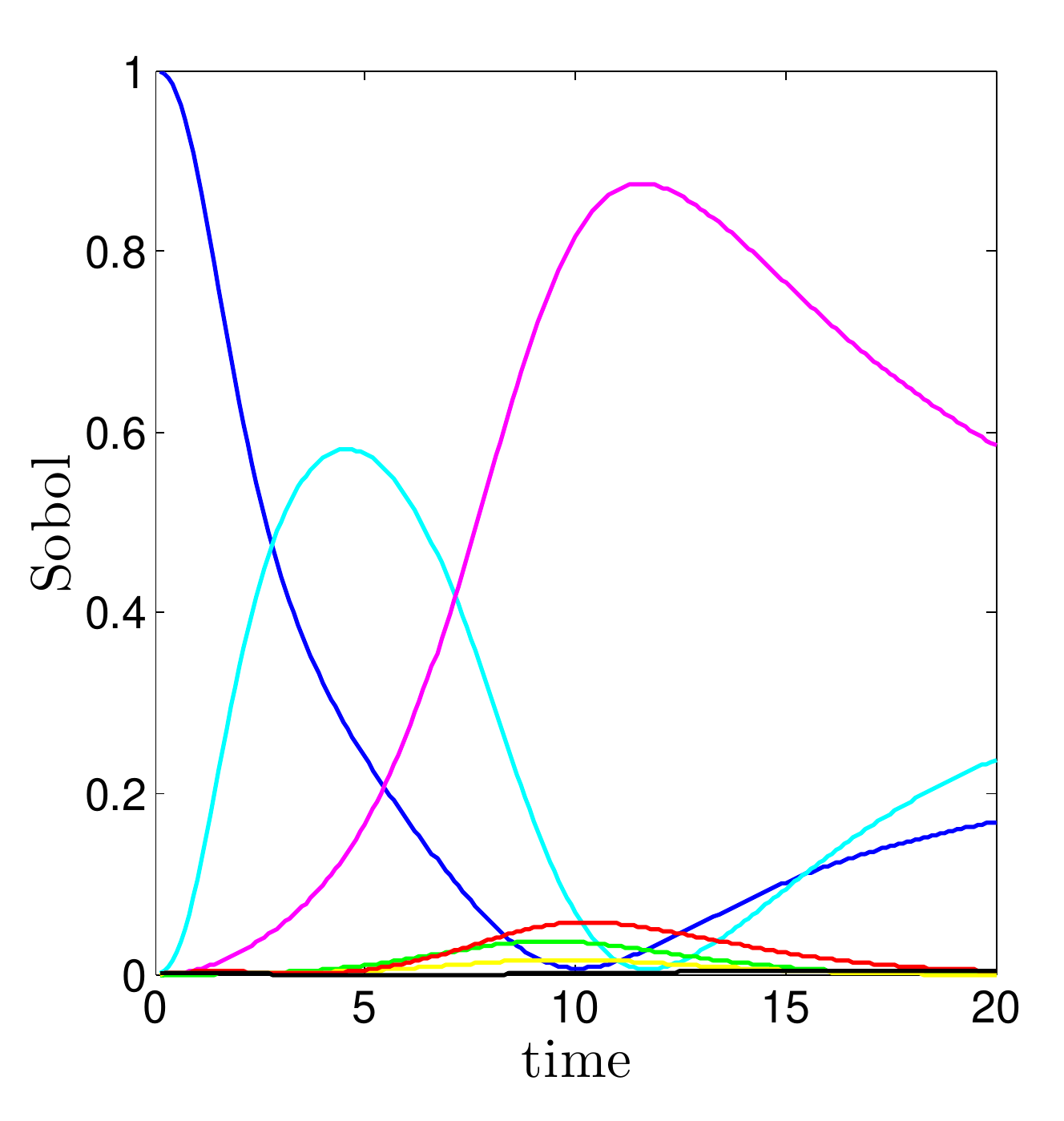}}
\subfloat[\c{P}]{
\includegraphics[width=0.33\linewidth]{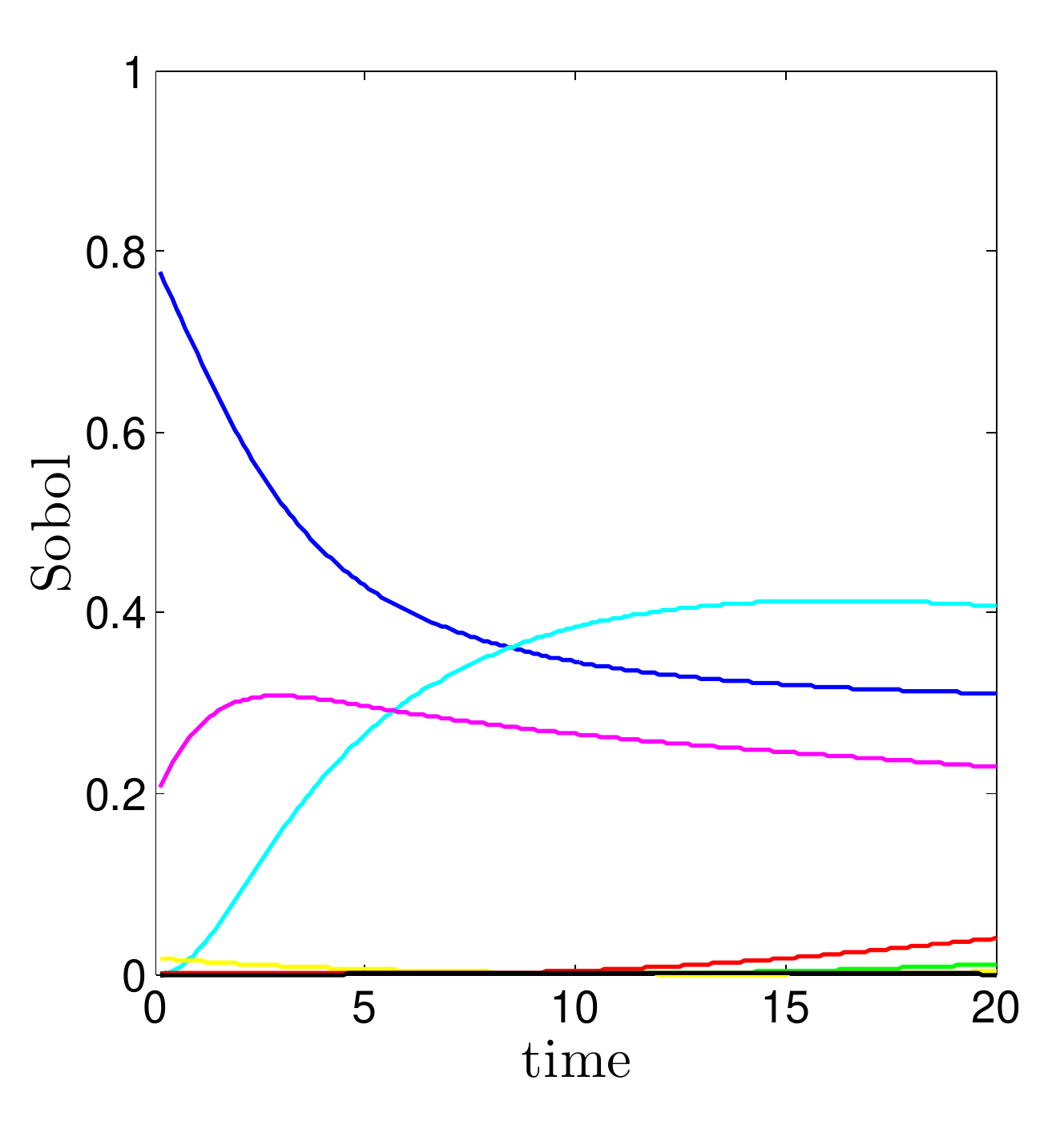}}
\end{tabular}
\caption{Dynamic plots of the Sobol' indices of the enzymatic reaction (\ref{equ:FEBS}) for the uncorrelated
case.}
\label{fig:sobolfebsU}
\end{figure}
\begin{figure}[htbp!]
\begin{tabular}{ccc}
\includegraphics[width=0.33\linewidth]{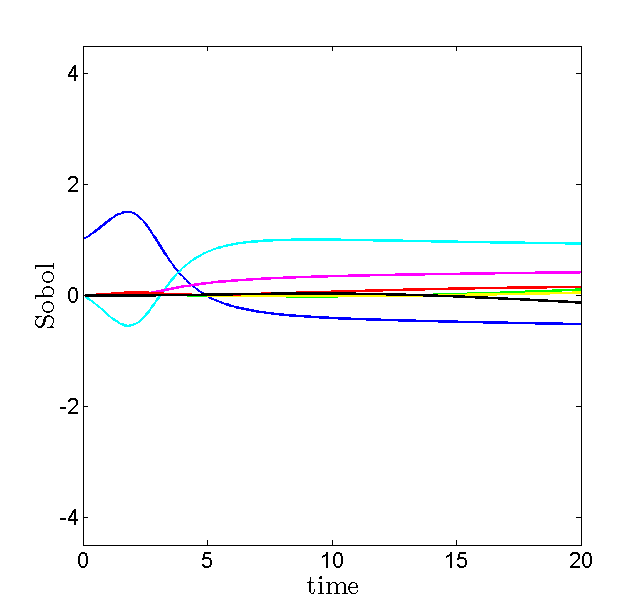}
\includegraphics[width=0.33\linewidth]{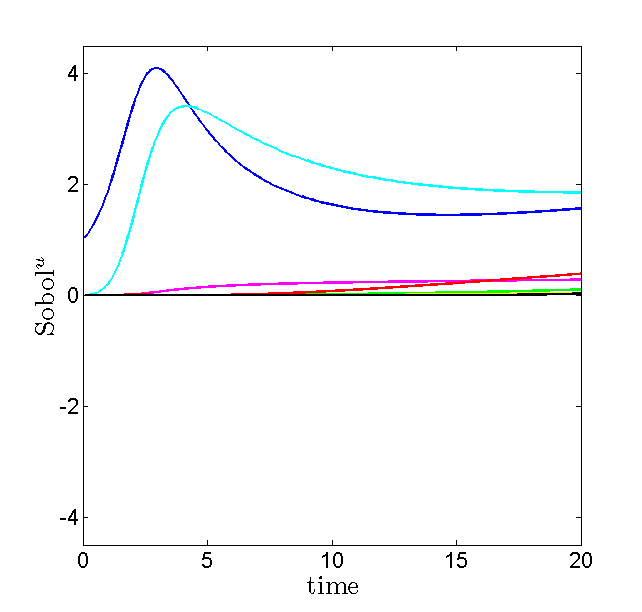}
\includegraphics[width=0.33\linewidth]{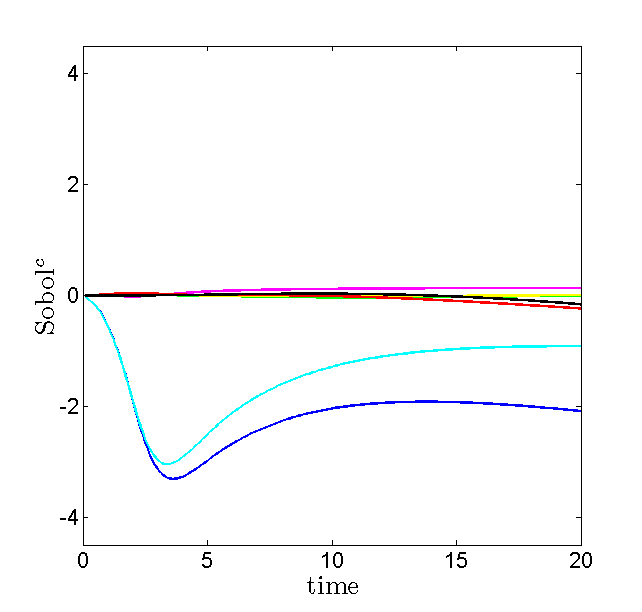}\\
\c{S}\\
\includegraphics[width=0.33\linewidth]{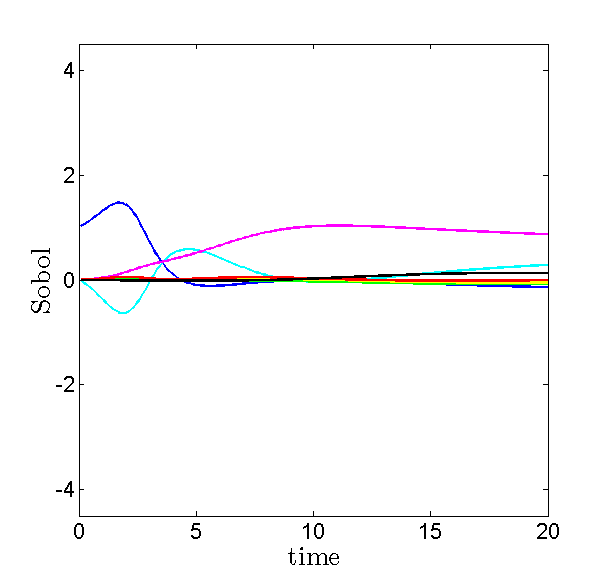}
\includegraphics[width=0.33\linewidth]{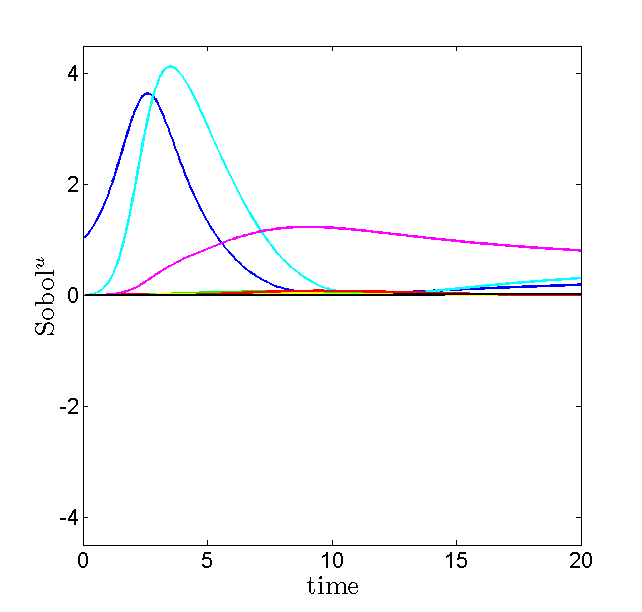}
\includegraphics[width=0.33\linewidth]{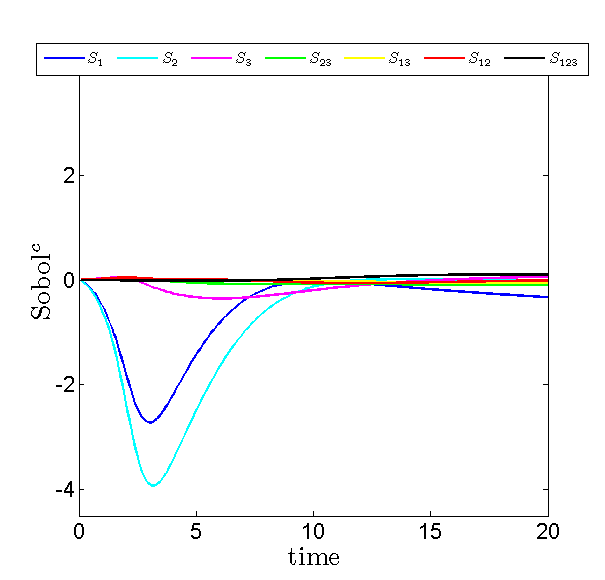}\\
\c{C} and \c{E}\\
\includegraphics[width=0.33\linewidth]{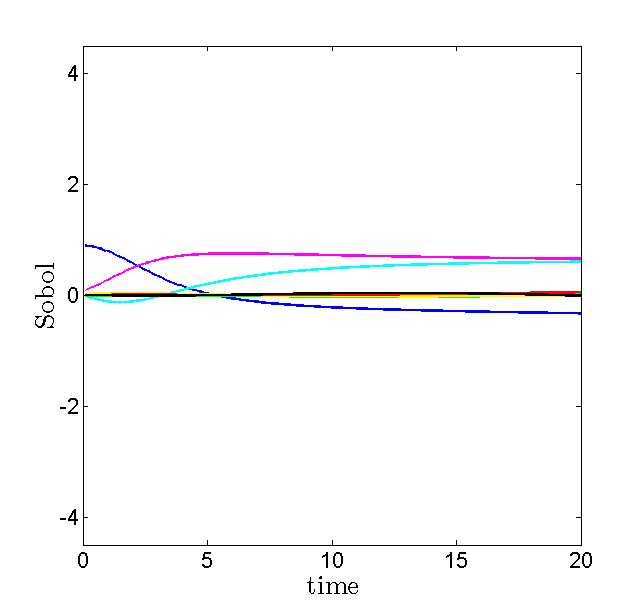}
\includegraphics[width=0.33\linewidth]{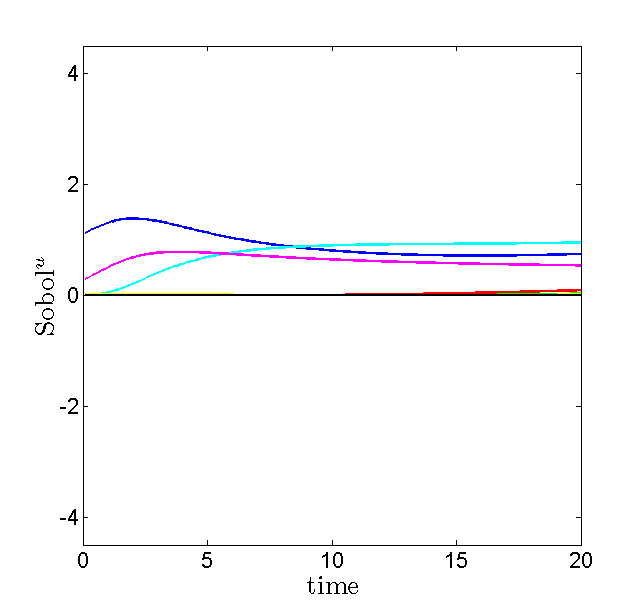}
\includegraphics[width=0.33\linewidth]{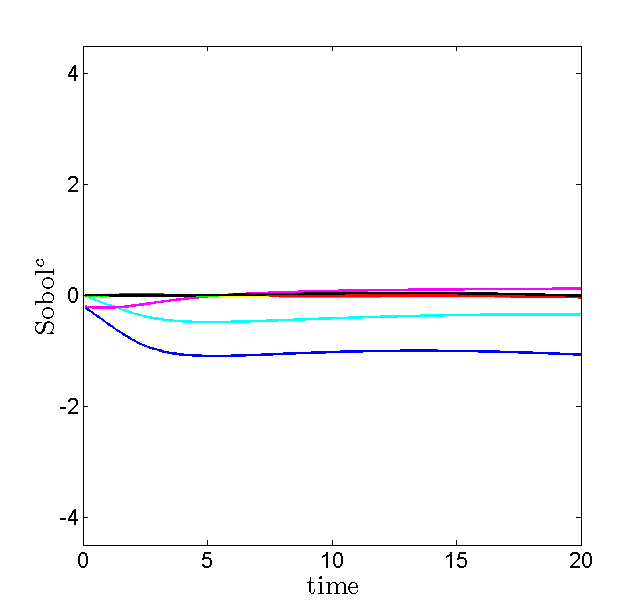}\\
\c{P}
\end{tabular}
\caption{Dynamic plots of the Sobol' indices of enzymatic reaction (\ref{equ:FEBS}), for the correlation 
matrix resulting from the experimental errors.
Left: Sobol' index $S$, middle: contribution from the respective variables, $S^u$, 
right: contribution due to correlation with other variables, $S^c$.}
\label{fig:sobolfebsC}
\end{figure}
Figures \ref{fig:Means}  and \ref{fig:Stds} show the obtained mean and standard deviation of all concentrations, 
respectively. Note, that the concentration of the complex \C\ has an equal dynamic 
behaviour - with opposite sign - as that of the enzyme \E. Again, it can be clearly seen from Fig.~\ref{fig:Stds}
that the amount of correlation between the random parameters has a significant impact on the variance of the
stochastic state variables. Without correlation the pre-equilibrium dynamics are seen only in the complex 
\C, and as a consequence in the enzyme \E, whereas the dynamics of the variance in the product \P\ simply
follow those of the substrate \S. With correlation the pre-equilibrium dynamics show up also in the substrate \S.
This difference is displayed too in the evolution of the Sobol' indices over time, see Figs.~\ref{fig:sobolfebsU}
and \ref{fig:sobolfebsC} (left column). From the plots in Fig. \ref{fig:sobolfebsC} it is again obvious 
that the interpretation of the Sobol' indices - and thus the ranking of the random variables - is less clear 
in case of correlated random variables. First of all, the Sobol' indices are no longer positive, and secondly, the 
contribution due to correlation can completely cancel the contribution from the variable itself, resulting
in a small Sobol' value, but fixing such a variable would have a large impact on the outcome.

\begin{table}[h!]
\centering
\begin{tabular}{|c|ccc|ccc|}
\hline
\c{S}& $S_i$& $S_i^\mathrm{u}$& $S_i^\mathrm{c}$& $S_i$&  $S_i^\mathrm{u}$& $S_i^\mathrm{c}$ \\
\hline
$S_1$&       0.35&  0.36& -0.01& -0.52&  1.57& -2.09\\
$S_2$&       0.45&  0.46& -0.01&  0.93&  1.85& -0.92\\
$S_3$&       0.08&  0.08& -0.00&  0.42&  0.28&  0.13\\
$S_{23}$&    0.02&  0.02& -0.00&  0.10&  0.10&  0.00\\
$S_{13}$&    0.00&  0.00& -0.00&  0.05&  0.04&  0.00\\
$S_{12}$&    0.09&  0.09& -0.00&  0.15&  0.39& -0.24\\
$S_{123}$&   0.00&  0.00&  0.00& -0.13&  0.03& -0.16\\
\hline
Sum&         0.99&  1.01& -0.02&  1.00&  4.26& -3.28\\  
\hline
\hline
\c{C}\ \& \c{E}& $S_i$& $S_i^\mathrm{u}$& $S_i^\mathrm{c}$& $S_i$&  $S_i^\mathrm{u}$& $S_i^\mathrm{c}$ \\
\hline
$S_1$&       0.17&  0.17& -0.00& -0.14&  0.19& -0.33\\
$S_2$&       0.24&  0.24& -0.00&  0.28&  0.31& -0.03\\
$S_3$&       0.58&  0.58& -0.00&  0.87&  0.80&  0.06\\
$S_{23}$&    0.00&  0.00& -0.00& -0.08&  0.01& -0.08\\
$S_{13}$&    0.00&  0.00&  0.00& -0.06&  0.00& -0.06\\ 
$S_{12}$&    0.00&  0.00&  0.00& -0.01&  0.00& -0.01\\
$S_{123}$&   0.00&  0.00&  0.00&  0.12&  0.02&  0.10\\
\hline
Sum&         0.99&  0.99&  0.00&  0.98&  1.33& -0.35\\   
\hline
\hline
\c{P}& $S_i$& $S_i^\mathrm{u}$& $S_i^\mathrm{c}$& $S_i$&  $S_i^\mathrm{u}$& $S_i^\mathrm{c}$ \\
\hline
$S_1$&       0.31&  0.31& -0.00& -0.33&  0.75& -1.08\\
$S_2$&       0.41&  0.41& -0.00&  0.61&  0.95& -0.34\\
$S_3$&       0.23&  0.23& -0.00&  0.66&  0.54&  0.12 \\
$S_{23}$&    0.01&  0.01& -0.00&  0.01&  0.03& -0.02\\
$S_{13}$&    0.00&  0.00&  0.00& -0.00&  0.01& -0.01\\
$S_{12}$&    0.04&  0.04& -0.00&  0.05&  0.09& -0.04 \\
$S_{123}$&  -0.00&  0.00& -0.00& -0.00&  0.00& -0.00\\
\hline
Sum&         1.00&  1.00& -0.00&  1.00&  2.37& -1.37\\  
\hline
\end{tabular}
\caption{Sobol' indices in the final time point $t_f=20$. The left part of the table corresponds 
to the uncorrelated case, correlation matrix $C_{u}$; 
the right part to the correlation matrix $C$ obtained from the FIM. 
The columns $S_{i}^{u}$ show  the contribution from the variables itself and
the columns $S_{i}^{c}$ the contribution from correlated variables; with exact computations the $S_i^c$-column 
in the left part of the table should be zero.}
\label{tb:Sobolfebs}
\end{table}
\begin{table}[h!]
\centering
\begin{tabular}{|c|ccc|ccc|}
\hline
\c{S}& $S_{T_{i}}$& $S_{T_{i}}^\mathrm{u}$& $S_{T_{i}}^\mathrm{c}$& $S_{T_{i}}$& $S_{T_{i}}^\mathrm{u}$& $S_{T_{i}}^\mathrm{c}$ \\
\hline
$S_{T_1}$&  0.44&   0.45& -0.01& -0.45&  2.03& -2.49\\
$S_{T_2}$&  0.56&   0.57& -0.01&  1.05&  2.37& -1.32\\
$S_{T_3}$&  0.10&   0.10&  0.00&  0.44&  0.45& -0.03\\                                    
\hline
\hline
\c{C}\ \& \c{E}& $S_{T_{i}}$& $S_{T_{i}}^\mathrm{u}$& $S_{T_{i}}^\mathrm{c}$& $S_{T_{i}}$& $S_{T_{i}}^\mathrm{u}$& $S_{T_{i}}^\mathrm{c}$ \\
\hline
$S_{T_1}$&  0.17&   0.17& -0.00& -0.09&  0.21& -0.30\\
$S_{T_2}$&  0.24&   0.24& -0.00&  0.31&  0.34& -0.02\\
$S_{T_3}$&  0.58&   0.58&  0.00&  0.85&  0.83&  0.02\\
\hline
\hline
\c{P}& $S_{T_{i}}$& $S_{T_{i}}^\mathrm{u}$& $S_{T_{i}}^\mathrm{c}$& $S_{T_{i}}$& $S_{T_{i}}^\mathrm{u}$& $S_{T_{i}}^\mathrm{c}$ \\
\hline
$S_{T_1}$&  0.35&   0.35& -0.00& -0.28&  0.85& -1.13\\
$S_{T_2}$&  0.46&   0.46& -0.00&  0.67&  1.07& -0.38\\
$S_{T_3}$&  0.24&   0.24&  0.00&  0.67&  0.58&  0.09 \\
\hline
\hline
\end{tabular}
\caption{Total Sobol' indices in the final time point $t_f=20$. 
The left part of the table corresponds to the uncorrelated case, correlation matrix $C_{u}$; 
the right part to the correlation matrix $C$ obtained from the FIM.}
\label{tb:SobolfebsT}
\end{table}
Ranking variables is often done based on the Sobol' indices. Depending on the time-frame of interest, one can
use as measure, e.g., the integral over time of the Sobol' indices, or the Sobol' indices in a specific point.
As an example of the latter we study at the final time-point, $t_f=20$, which part of the variance of the 
concentrations is due to the variance of respective input variables. 
Tables~\ref{tb:Sobolfebs} and \ref{tb:SobolfebsT} give the Sobol' indices
and the total Sobol' indices, respectively, for the random variables $k_1$, $k_2$, and $k_3$. Note that the
sum of the Sobol' indices equals one and for the uncorrelated case the $S^c$ values should be zero, so the 
deviation of these values gives the error in the approximation, in this case mostly due to the approximation of
the high-order moments needed for the Sobol' indices. From Table~\ref{tb:SobolfebsT} one can deduce that, in
the uncorrelated case, decomposing 
the variance of the concentration of the substrate \S\  would lead to the conclusion that
the random variable $k_3$ can be fixed to the nominal deterministic value, since its influence on the
variance of \S\ is negligible. In the correlated case, however, this is no longer true. One can also see
that especially for $k_1$ and $k_2$ the interpretation of the Sobol' indices is not trivial, e.g., for \C\ the
influence of $k_1$ is very small but this is due to cancellation of the $S^u$ and $S^c$ contribution.
As stated before, the interpretation of the Sobol' indices for correlated random variables is largely an open
issue, but it is clear that in general less random input variables can be fixed to their nominal deterministic
values.

\section{Discussion and concluding remarks}
In this paper we extended the PCE method to general multivariate distributions, including correlations
between the random input variables, by constructing an orthogonal polynomial basis for the probability
space. We showed that the usual statistics like mean, variance and Sobol' indices can be obtained with
no extra cost, as is the case for the current PCE implementations. This method makes it possible to compute
for a set of random input variables with any multivariate distribution, including an experimentally 
determined one, the stochastic distribution of a Quantity of Interest. An application we studied is 
the true propagation of experimental errors through the parameter-fitting process onto the state variables,
so that a QoI based on these state variables is described as a distribution and, e.g.,
optimal experimental design can be applied to reduce the variance in such a QoI.
In this paper we used either exact integration or Monte Carlo to compute the moments and the projection onto
the polynomial basis, but for the latter Gauss quadrature is under development.

\medskip
There are a few remarks to make with respect to the method
\begin{itemize}
 \item As for other PCE methods, the polynomial base is only dependent on the multivariate distribution
 of the random input variables. Once computed, all problems with random input variables described by
 this pdf can be solved using the same base. 
 \item To obtain the polynomial basis raw moments are needed which have to be computed accurately enough 
 not to pollute the results. This is not a trivial task for high dimensional problems or a high polynomial order. E.g., in our first example, the size of 
 the moments range in order of magnitude from 1 for the first ten orders to $10^{15}$ for order 75 for 
 an uncorrelated Normal distribution. If the correlation coefficient $\varrho$ equals $0.9$ the latter 
 order of magnitude is even much higher: $10^{21}$ for moments of order 75.
 \item The orthogonal polynomial basis is not unique, it is dependent on the choice and ordering of
 the original set of linearly independent polynomials. It is an open question at the moment whether there
 is an optimal choice for this original set of polynomials, e.g., to reduce the magnitude of the moments 
 required or to simplify computations.
% \item Different basis for G-S to reduce size of moments? (uncorrelated is much smaller than correlated)
% \item Gauss quadrature based on new polynomials
 \item Whereas in the uncorrelated case interpretation of Sobol' indices and ranking of the influence
 of the variation of the random input variables on the variance in the QoI is more or less straighforward,
 this is no longer true for correlated random variables. 
\end{itemize}

\section*{Acknowledgements}
MN and JB acknowledge support by European Commission's 7th Framework Program, project BioPreDyn grant number 289434.

 \end{document}